\documentclass[12pt]{article}
\usepackage{amsfonts,amsmath,amssymb,mathrsfs} 
\usepackage{mathtools}
\usepackage{latexsym}
\usepackage{mleftright}
\usepackage{cite}
\usepackage{txfonts}
\allowdisplaybreaks
\newtheorem{Theorem}{Theorem}[section]
\newtheorem{Lemma}{Lemma}[section]
\newtheorem{Remark}{Remark}[section]

\newcommand{\beq}{\begin{eqnarray}}
\newcommand{\eeq}{\end{eqnarray}}
\newcommand{\beqno}{\begin{eqnarray*}}
\newcommand{\eeqno}{\end{eqnarray*}}
\newcommand{\be}{\begin{equation}}
\newcommand{\ee}{\end{equation}}
\newcommand{\ideq}{{\lower .5ex
\hbox{$\>\>\stackrel{\triangle}{=}\>\>$} }}

\def \d2{\Delta_{+}^2}

\begin{document}
\newcommand{\D}{\displaystyle}
\title{\bf Strong Solutions for 1D Compressible Navier-Stokes/Allen-Cahn System with Phase Variable Dependent Viscosity}
\author{ Yuanxiang Yan, Shijin Ding and Yinghua Li\thanks{Corresponding author.
Email: yuanxiangyan@m.scnu.edu.cn (Y. Yan), dingsj@scnu.edu.cn (S. Ding), yinghua@scnu.edu.cn (Y. Li).}
\\
\
\\
{\footnotesize\it School of Mathematical Sciences,
South China Normal University, Guangzhou 510631, China}
}
\date{}
\maketitle

\allowdisplaybreaks

\begin{abstract}
This paper is concerned with a non-isentropic compressible Navier-Stokes/Allen-Cahn system
with phase variable dependent viscosity $\eta(\chi)=\chi^\alpha$
and temperature dependent heat-conductivity $\kappa(\theta)=\theta^\beta$.
We show the global existence of strong solutions under some assumptions on growth exponent $\alpha$ and initial data.
It is worth noting that the initial data could be large if $\alpha\ge0$ is small,
and the growth exponent $\beta>0$ can be arbitrary large.
\end{abstract}

existence; global solutions.

\section{Introduction}
In this paper, we investigate a diffuse interface model for two-phase flows of viscous compressible fluids,
which was proposed by Blesgen\cite{TB}.
This model can be used to describe topological transitions on the interface such as droplet coalescence or droplet break-up.
A lot of attentions have been paid to diffuse interface models because of their clear background
and their applications in numerical simulations.
Great progresses have been achieved in the studies of incompressible case, i.e. $\rho=const.$, see \cite{AH09, G-G10, G-M-T} for example.
The researches on compressible diffuse interface models mainly focus on
Navier-Stokes/Allen-Cahn system \cite{TB, H-M-R} and Navier-Stokes/Cahn-Hilliard system \cite{A-F, L-T}.
The theoretical analysis of compressible Navier-Stokes/Allen-Cahn model began with
Feireisl et al. \cite{F-P-R-S} and Kotschote \cite{Kotschote_1}.
They proved the existence of weak solutions for isentropic system and local strong solutions for non-isentropic system, respectively.
Here, we are interested in non-isentropic compressible Navier-Stokes/Allen-Cahn system,
which was simplified by Chen et al. \cite{C-H-H-S-1,C-H-H-S-2} into the following form
\begin{align*}
\begin{cases}
\rho_t+{\rm div}(\rho{\bf u})=0,
\\
\rho{\bf u}_t+\rho({\bf u}\cdot\nabla){\bf u}-2{\eta(\chi)}{\rm div}{\mathbb D}{\bf u}
-{\lambda(\chi)}\nabla{\rm div}{\bf u}
=-{\rm div}\Big(\delta\nabla\chi\otimes\nabla\chi-\dfrac\delta2|\nabla\chi|^2+\theta\dfrac{\partial p}{\partial\theta}\Big),
\\
\rho\chi_t+\rho({\bf u}\cdot\nabla)\chi=-\mu,
\\
\rho\mu=\dfrac\rho\delta(\chi^3-\chi)-\delta\Delta\chi,
\\
c_v(\rho\theta_t+\rho{\bf u}\cdot\nabla\theta)+\theta p_\theta{\rm div}{\bf u}-{\rm div}(\kappa(\theta)\nabla\theta)
=2{\eta(\chi)}|{\mathbb D}{\bf u}|^2+{\lambda(\chi)}({\rm div}{\bf u})^2+\mu^2,
\end{cases}
\end{align*}
where $\rho, {\bf u}, \chi, \theta$ represent the total density, the mean velocity of the fluid mixture, the phase field variable.
Moreover, $\mu$ is the chemical potential, $\delta$ and $c_v$ are related to the thickness of
the interfacial region and the heat capacity at constant volume.
The viscosity coefficients $\eta(\chi)$ and $\lambda(\chi)$ satisfy
$\eta>0, \lambda + \dfrac{2}{N}\eta>0$.
$\kappa(\theta)$ is the heat-conductivity, and the pressure $p=R\rho \theta$ with $R>0$.

\vskip2mm
In this paper, we assume that the viscosity coefficients satisfy $\lambda(\chi)=-\eta(\chi)$, then the above system in 1D becomes
\begin{equation}
\label{NSAC1D}
\begin{cases}
\rho_t+(\rho u)_{\tilde{x}}=0,
\\
\rho u_{t}+\rho uu_{\tilde{x}}+(R\rho\theta)_{\tilde{x}}=(\eta(\chi) u_{\tilde{x}})_{\tilde{x}}-\D\frac\delta2\left(\chi_{\tilde{x}}^{2}\right)_{\tilde{x}},
\\
\rho\chi_{t}+\rho u\chi_{\tilde{x}}=-\mu,
\\
\D\rho\mu=-\delta\chi_{{\tilde{x}}{\tilde{x}}}+\frac\rho\delta(\chi^3-\chi),
\\
c_v(\rho \theta_t + \rho u \theta_{\tilde{x}})+ R\rho\theta u_{\tilde{x}} - (\kappa(\theta)\theta_{\tilde{x}})_{\tilde{x}} = \eta u_{\tilde{x}}^2 + \mu^2
\end{cases}
\end{equation}
for $({\tilde{x}},t)\in(0,1)\times(0, +\infty)$.
We supplement (\ref{NSAC1D}) with initial value conditions
\begin{gather}
\label{I}
(\rho,\, u, \, \chi, \,\theta)\Big|_{t=0}=(\rho_0,\, u_0, \, \chi_0, \,\theta_0), \qquad {\tilde{x}}\in(0, 1)
\end{gather}
and boundary value conditions
\begin{gather}
\label{B1}
(u, \, \chi_{\tilde{x}},\,\theta_{\tilde{x}})\Big|_{{\tilde{x}}=0,1}=(0,\,0, \, 0), \qquad  t\ge0.
\end{gather}
Without loss of generality, we assume
$\displaystyle\int_0^1\rho_0({\tilde{x}})d{\tilde{x}}=1$.
Then in Lagrange coordinates
$x = \displaystyle\int_0^{\tilde{x}}\rho (\xi,t)d{\xi}$,
the system $(\ref{NSAC1D})$--$(\ref{B1})$ can be rewritten as
\begin{equation}
\label{L-NSAC1D}
\begin{cases}
v_t = u_x,
\\[0.5em]
u_{t}+\left(\dfrac{\theta}{v}\right)_{x}=\left(\dfrac{\eta(\chi) u_{x}}{v}\right)_{x}-\dfrac{1}{2}\left(\dfrac{\chi_{x}^2}{v^2}\right)_{x},
\\[0.5em]
\chi_{t}=-v\mu,
\\[0.5em]
\D\mu=-\left(\frac{\chi_{x}}{v}\right)_x+ (\chi^3-\chi),
\\[0.8em]
\theta_t + \dfrac{\theta}{v} u_{x}  = \left(\dfrac{\kappa(\theta)\theta_{x}}{v} \right)_{x} + \dfrac{\eta u_{x}^2}{v} + v\mu^2
\end{cases}
\end{equation}
with initial and boundary value conditions
\begin{align}
\label{IC}
&(v,\, u, \, \chi, \,\theta)\Big|_{t=0}=(v_0,\, u_0, \, \chi_0, \,\theta_0), \quad {x}\in(0, 1),
\\
\label{BC}
&(u, \, \chi_{x},\,\theta_{x})\Big|_{x=0,1}=(0,\,0, \, 0), \qquad~~~~~  t\ge0.
\end{align}
Here $v=\dfrac{1}{\rho}$ represents specific volume.
In the following, we choose
$$
\eta(\chi)=\tilde{\eta}\,\chi^{\alpha},\qquad \kappa(\theta)=\tilde{\kappa}\,\theta^\beta,
$$
with $\alpha \ge 0,~ \beta>0$, and the constants $R = c_v=\tilde{\eta} = \tilde{\kappa} = \delta = 1$.

\vskip2mm

Before introducing our main result, we first give a brief review on some related works.
For isentropic compressible NSAC model with constant viscosity,
Ding et al. \cite{D-L-L} and Chen et al. \cite{C-G} obtained the 1D global well-posedness
without vacuum and with vacuum, respectively.
Later, Ding et al. \cite{D-L-T} derived the existence and uniqueness of global strong solutions with free boundary condition.
Very recently, Chen and Zhu \cite{C-Z} assumed that the viscosity coefficient satisfied
$$
\eta(\rho,\chi)=1+\rho^\alpha\chi^\beta.
$$
They proved the existence and uniqueness of global classical solutions when $2\le \alpha\le \gamma$ and $\beta=0$.
In the case of $\beta\ge1$, they obtained a blow-up criterion for strong solutions.
In \cite{C-G, C-Z}, the phase variable $\chi$ was assumed to satisfy Dirichlet boundary conditions, i.e. $\chi\big|_{x=0,1}=0$.
For non-isentropic compressible NSAC system with constant viscosity,
Chen et al. \cite{C-H-H-S-1,C-H-H-S-2} studied global strong solutions of initial-boundary value problem and Cauchy problem.
Besides, there are also some researches on asymptotic behavior, weak solutions, stationary solutions and wave problem,
one can find them in the references \cite{C-H-S, C-W-Z19, C-W-Z20, Kotschote_2, L-Y-Z, L-Y-Z20, Y-Z}.

\vskip2mm
The main purpose of this paper is to deal with phase variable dependent viscosity.
Since this model is governed by the full compressible Navier-Stokes equation coupled with Allen-Cahn equations,
let's do some review on the full compressible Navier-Stokes equations.
When both the viscosity $\eta$ and the heat-conductivity $\kappa$ are positive constants,
the analysis mainly relies on the upper and lower bounds of specific volume $v$ and temperature $\theta$.
The proof was built upon a representation of specific volume $v$, which was obtained by Kazhikhov et al. \cite{K1,K2}.
When the viscosity $\eta$ is a constant or depends only on $v$, mass conservation equation and momentum conservation equation imply
\begin{equation}\label{eta-v}
\left( \dfrac{\eta(v)v_x}{v} \right)_t = u_t + P_x,
\end{equation}
which was observed by Kanel in \cite{Kanel}. By virtue of $(\ref{eta-v})$, one can get global well-posedness of solutions with large initial data,
see \cite{C-H-T,Daf,D-H,Kawohl,Q-Y,T-Y-Z-Z,J} and references therein.
When the viscosity $\eta$ depends on the temperature $\theta$ and the specific volume $v$,
the identity $(\ref{eta-v})$ becomes
\begin{equation}\label{eta-v theta}
\left( \dfrac{\eta(v,\theta)v_x}{v} \right)_t = u_t + P_x + \dfrac{\eta_\theta(v,\theta)}{v}\left(\theta_t v_x - u_x \theta_x \right).
\end{equation}
It is clear that the temperature dependent viscosity has a strong influence on the solution.
And as pointed out in \cite{J-K}, such a dependence turns out to be challenging.
Later, Wang and Zhao \cite{W-Z} considered the following case
\begin{equation*}
\eta(v,\theta)= \tilde{\eta}h(v)\,\theta^\alpha, \qquad  \kappa(\theta)=\tilde{\kappa}h(v)\,\theta^\alpha
\end{equation*}
under some structure assumptions. But their results excluded the case of $h\equiv{const.}$
Recently, Sun, Zhang and Zhao \cite{S-Z-Z} assumed that
\begin{equation*}
\eta(\theta)=\theta^{\alpha},\qquad \kappa(\theta)=\theta^\beta,
\end{equation*}
with $\alpha\ge0, \beta \ge0$ and the initial data $v_0 \ge V_0, \theta_0 \ge V_0$ for some constant $V_0>0$.
When $\alpha$ is small, they obtained the existence and uniqueness of strong solutions.
Thanks to the ideas in \cite{S-Z-Z}, we can handle our problem for the case $\eta(\chi)=\chi^\alpha$.

\vskip4mm
Our main result is the following global-in-time existent theorem.
\begin{Theorem}
\label{th:1.1}
For given positive numbers $M_0$ and $V_0$, assume that
\begin{equation}\label{local C}
 v_0 \ge V_0,\quad
 V_0\le\chi_0 \le 1,\quad
 \theta_0 \ge V_0, \quad
 ||(v_0,u_0,\theta_0)||_{H^2}+ ||\chi_0||_{H^3} \le M_0.
\end{equation}
Then there exists a positive constant $\epsilon_0>0$, depending only on $M_0$, $V_0$ and $\beta$, such that the problem $(\ref{L-NSAC1D})$--$(\ref{BC})$ with $0\le \alpha \le \epsilon_0$ and $\beta>0$ admits a unique global strong solution $(v,u,\chi,\theta)$ on $({x},t)\in[0,1]\times[0, +\infty)$, satisfying
\begin{equation*}
\inf\limits_{({x},t)\in[0,1]\times[0, +\infty)}\{v(x,t),\theta(x,t)\}>0, \quad
\sup\limits_{({x},t)\in[0,1]\times[0, +\infty)}\{v(x,t),\theta(x,t)\}<\infty,
\end{equation*}
\begin{align*}
V_0\le\chi(x, t)\le 1,\qquad (x, t)\in[0,1]\times[0, +\infty),
\end{align*}
and
\begin{align*}
&(v,u,\theta)\in C([0,\infty);H^2),\quad \chi \in C([0,\infty);H^3),\\
&v_x\in L^2(0,\infty;H^1),\quad (u_x, \chi_x, \theta_x)\in L^2(0,\infty;H^2).
\end{align*}
\end{Theorem}

\begin{Remark}
$(i)$ Even though we can just handle the case when $\alpha$ is small and require $\chi$ has positive lower bound,
but this is a first result on global strong solutions with phase variable dependent viscosity,
which is more in line with physical reality.
\\
$(ii)$ The strong solutions exist on whole time interval $[0, +\infty)$, so all the bounds in a priori estimates are time-independent.
This makes it possible to consider long time behavior of the solutions.
\\
$(iii)$ The theorem is also correct if one replaces $\eta=\eta(\chi)=\chi^\alpha$ by $\eta=\eta(\theta)=\theta^\alpha$.
\end{Remark}

The key step of the proof is to obtain lower and upper bounds of the specific volume $v$.
We apply the argument developed by Kazhikhov \cite{K1} to derive a presentation of specific volume (Lemma $\ref{Lemma-v}$).
In our model, $(\ref{eta-v theta})$ becomes
\begin{equation}\label{eta chichi}
\left( \dfrac{\eta(\chi)v_x}{v} \right)_t = u_t + P_x + \dfrac{1}{2}\left( \dfrac{\chi_x^2}{v^2} \right)_x+\dfrac{\alpha \chi^{\alpha-1}}{v}\left(\chi_t v_x - u_x \chi_x \right).
\end{equation}
The last term in ($\ref{eta chichi}$) is highly nonlinear.
Fortunately, we can use the smallness of $\alpha$ to control it as in \cite{S-Z-Z}.
Moreover, comparing with ($\ref{eta-v theta}$), we need to deal with the high-order strongly nonlinear term
$\left( \dfrac{\chi_x^2}{v^2} \right)_x$ (Lemma \ref{L-chit}, Lemma \ref{L-chixx}).

\vskip2mm
The structure of this paper is as follows.
In Section 2, we do some a priori estimates which are independent with $T$.
Then We finish the proof of Theorem $\ref{th:1.1}$ by using continuity method in Section $\ref{section proof}$.


\setcounter{section}{1}
\setcounter{equation}{0}
\section{A priori estimates}\label{section priori}

\allowdisplaybreaks
For given some positive constants $m_i(i=1,2,3)$ and $N$, we define the set
\ \par
\[
X(0,T;m_1,m_2,m_3,N) := \mleft\{\,(v,u,\chi, \theta) \,\middle|
\begin{array}{c}
(v,u,\theta)\in C([0,T]; H^2),~\chi \in C([0,T]; H^3),\\
(v_t, \chi_t) \in L^\infty(0,T; H^1), \\
(v_x,u_t,\theta_t) \in L^2(0,T; H^1),\\
(u_x,\theta_x,\chi_x,\chi_t)\in L^2(0,T; H^2),\\
\mathcal{E}(0,T)\le N^2,\\
v\ge m_1,\chi\ge m_2,\theta\ge m_3,
\forall (x,t)\in [0,1]\times[0,T].
\end{array}
\mright\}
\]
where
$$
\mathcal{E}(0,T):= \sup\limits_{0\le t \le T}||(v_x,u_x,\chi_x, \theta_x)(t)||_{H^1}^2
+ \int_0^T ||\chi_t||_{L^2}^2 dt.
$$

Our main purpose of this section is to derive the time-independent a priori estimates of the solutions $(v,u,\theta,\chi)\in X(0,T;m_1,m_2,m_3,N)$ to the problem $(\ref{L-NSAC1D})$--$(\ref{BC})$ for $0<m_i\le 1~(i=1,2,3),~N\ge 8$ and $T>0$. From here to the end of this paper, $C$ and $C_i~(i=1,2,\cdots)$ denote the generic positive constants, dependenting only on $\beta$, $V_0$ and $M_0$.

In order to simplify the presentation, without loss of generality, we assume that
\begin{equation}\label{EE}
\int_0^1 v_0 dx = 1,\quad \int_0^1 \left( \theta_0 + \dfrac{u_0^2}{2}+ \dfrac{(\chi_0^2-1)^2}{4}+\dfrac{\chi_{0x}^2}{2v_0} \right) dx =1.
\end{equation}

\begin{Lemma}
\label{L-energy}
Let $(v,u,\chi,\theta)\in X(0,T;m_1,m_2,m_3,N)$ be a solution to the problem $(\ref{L-NSAC1D})$--$(\ref{BC})$ on $[0,1]\times[0,T]$. Then
\begin{equation}
\label{energy}
\sup\limits_{0\le t \le T}\int_0^1\left( \dfrac{u^{2}}{2}+\dfrac{(\chi^2-1)^2}{4}+\dfrac{\chi_{x}^2}{2v}+\Phi(v)+\Phi(\theta)\right)dx+\int_{0}^{T}W(\tau) d\tau
\le E_0,
\end{equation}
where
\begin{equation*}
\Phi(s)= s-\ln s -1,\quad  W(t) = \int_0^1 \left( \dfrac{\theta^\beta\theta_x^2}{v\theta^2}+\dfrac{\eta(\chi)u_x^2}{v\theta}+\dfrac{v\mu^2}{\theta} \right)dx,
\end{equation*}
and
\begin{equation*}
E_0=\int_0^1\left(\dfrac{u_{0}^{2}}{2}+ +\dfrac{\left(\chi_{0}^{2}-1\right)^{2}}{4}+\dfrac{\chi_{0x}^{2}}{2v_0}+ \Phi(v_0)+\Phi(\theta_0)\right)dx.
\end{equation*}
\end{Lemma}
\noindent{\it\bfseries Proof.}\quad
From $(\ref{BC})$, multiplying $(\ref{L-NSAC1D})_3$ by $\mu$, adding $(\ref{L-NSAC1D})_1$ and $(\ref{L-NSAC1D})_5$ integrating them over $[0,1]$ by parts, combining with $(\ref{EE})$, we have
\begin{equation}\label{v=1}
\int_0^1 v\,dx = 1,\quad \int_0^1\left(\theta + \dfrac{u^{2}}{2}+\dfrac{(\chi^2-1)^2}{4}+\dfrac{\chi_{x}^2}{2v}\right)dx =1.
\end{equation}
Multiplying $(\ref{L-NSAC1D})_1$, $(\ref{L-NSAC1D})_2$, $(\ref{L-NSAC1D})_3$ and $(\ref{L-NSAC1D})_5$ by $1-v^{-1}$, $u$, $\mu$ and $1-\theta^{-1}$, respectively, integrating by parts over $[0,1]$, adding them together and combining the boundary condition $(\ref{BC})$, we get
\begin{equation*}
\frac{d}{dt}\int_0^1\left( \dfrac{u^{2}}{2}+\dfrac{(\chi^2-1)^2}{4}+\dfrac{\chi_{x}^2}{2v}+(v-\ln v)+(\theta-\ln \theta)\right)dx+ W(t) = 0.
\end{equation*}
then the proof of Lemma $\ref{L-energy}$ is finished.
\hfill$\Box$

\begin{Lemma}\label{Lemma-v}
Let $(v,u,\chi,\theta)\in X(0,T;m_1,m_2,m_3,N)$ be a solution to the problem $(\ref{L-NSAC1D})$--$(\ref{BC})$ on $[0,1]\times[0,T]$, we define $\eta_0:=\eta(\chi_0)$, then for any $t\ge0$, there is a $\alpha_0(t)\in (0,1)$ such that
\begin{equation}\label{Lv}
v = B(t)D(x,t) + \int_0^t \dfrac{B(t)D(x,t)}{B(\tau)D(x,\tau)}v(x,\tau)J(x,\tau)\,d\tau,
\end{equation}
where
\begin{align}
\label{g}
&g(x,t):=-\left[ u \left( \dfrac{1}{\eta}\right)_{t}+ \dfrac{\theta}{v}  \left( \dfrac{1}{\eta}\right)_{x} + \dfrac{\chi_{x}^2}{2v^2} \left( \dfrac{1}{\eta}\right)_{x} + \dfrac{\eta_x u_x}{\eta v} \right],
\\
\label{B}
&B(t):=\exp\left\{ -\int_0^t \int_0^1 \left( \dfrac{\theta+u^2}{\eta}+ \dfrac{\chi_x^2}{2\eta v}\right)dx d\tau \right\},
\\
\label{D}
&D(x,t):=v_0(x) \exp\left\{ \int_{\alpha_0(t)}^x \dfrac{u}{\eta}dy - \int_{0}^x \dfrac{u_0}{\eta_0}dy + \int_0^1 v_0 \left(  \int_{0}^x \dfrac{u_0}{\eta_0}dy \right)dx \right\},
\\
\label{J}
&J(x,t):=\left( \dfrac{\theta}{\eta v}+ \dfrac{\chi_x^2}{2\eta v^2} + \int_0^x g dy \right)(x,t) - \int_0^1 v \left( \int_0^x g dy \right)(x,t)dx .
\end{align}
\end{Lemma}

\noindent{\it\bfseries Proof.}\quad
For $g=g(x,t)$ given by  $(\ref{g})$, it follows from $(\ref{L-NSAC1D})_1$ and $(\ref{L-NSAC1D})_2$ that
\begin{equation}\label{lnv_xt}
\left( \dfrac{u}{\eta} \right)_{t} + \left( \dfrac{\theta}{\eta v} \right)_{x} + \left( \dfrac{\chi_x^2}{2 \eta v^2} \right)_{x} + g(x,t) = \left( \dfrac{u_x}{v} \right)_{x}=(\ln v)_{xt}.
\end{equation}

Define
\begin{equation}\label{varphi}
\varphi(x,t):= \int_0^t \left( \dfrac{u_x}{v} - \dfrac{\theta}{\eta v} - \dfrac{\chi_x^2}{2 \eta v^2} - \int_0^x g(y,\tau) dy \right) (x,\tau) d\tau + \int_0^x  \left( \dfrac{u_0}{\eta_0}\right)(y)dy,
\end{equation}
where $\eta_0:=\eta(\chi_0)$. Then we have
\begin{equation*}
\varphi_t = \dfrac{u_x}{v} - \dfrac{\theta}{\eta v} - \dfrac{\chi_x^2}{2 \eta v^2} - \int_0^x g(y,t) dy,
\qquad  \varphi_x = \dfrac{u}{\eta}.
\end{equation*}
Combining with $(\ref{L-NSAC1D})_1$, we arrive at
\begin{equation}\label{vuvarphi}
(v\varphi)_{t} - (u\varphi)_x = u_x - \dfrac{\theta+u^2}{\eta}- \dfrac{\chi_x^2}{2\eta v}- v\int_0^x g(y,t) dy.
\end{equation}
Integrating $(\ref{vuvarphi})$ over $[0,1]\times[0,t]$ and using $(\ref{BC})$, we get
\begin{equation*}
\int_0^1 (v\varphi)dx -\int_0^1 (v_0\varphi_0)dx = - \int_0^t\int_0^1 \left( \dfrac{\theta+u^2}{\eta}+ \dfrac{\chi_x^2}{2\eta v}+ v\int_0^x g(y,t) dy \right)dxd\tau.
\end{equation*}

On one hand, thanks to the mean value theorem and $(\ref{v=1})$, there is a $\alpha_0(t)$ for any $t\ge0$ such that
\begin{align} \label{varphi_1}
&\varphi(\alpha_0(t),t)= \int_0^1 (v\varphi)dx \\
=&\int_0^1 v_0\left(\int_0^x \dfrac{u_0}{\eta_0}dy \right)dx - \int_0^t\int_0^1 \left( \dfrac{\theta+u^2}{\eta}+ \dfrac{\chi_x^2}{2\eta v}+ v\int_0^x g(y,t) dy \right)dxd\tau. \nonumber
\end{align}
On the other hand, it follows from $(\ref{varphi})$ that
\begin{equation} \label{varphi_2}
\begin{split}
&\varphi(\alpha_0(t),t) = \ln v(\alpha_0(t),t)-\ln v_0(\alpha_0(t))
-\int_0^t \left( \dfrac{\theta}{\eta v}\right)(\alpha_0(t),\tau) d\tau 
\\
& \qquad-\int_0^t \left(\dfrac{\chi_x^2}{2 \eta v^2}\right)(\alpha_0(t),\tau) d\tau
 -\int_0^t\int_0^{\alpha_0(t)} g(y,\tau) dy d\tau
+\int_0^{{\alpha_0(t)}} \left( \dfrac{u_0}{\eta_0}\right)(y)dy.
\end{split}
\end{equation}
Hence, collecting $(\ref{varphi_1})$ and $(\ref{varphi_2})$, we have
\begin{equation} \label{lnv-lnv_0}
\begin{split}
&\ln v(\alpha_0(t),t)-\ln v_0(\alpha_0(t))
\\
=&\int_0^1 v_0\left(\int_0^x \dfrac{u_0}{\eta_0}dy \right)dx-\int_0^t\int_0^1 \left( \dfrac{\theta+u^2}{\eta}+ \dfrac{\chi_x^2}{2\eta v}+ v\int_0^x g(y,t) dy \right)dxd\tau
\\
&+\int_0^t \left( \dfrac{\theta}{\eta v}\right)(\alpha_0(t),\tau) d\tau
+\int_0^t \left(\dfrac{\chi_x^2}{2 \eta v^2}\right)(\alpha_0(t),\tau) d\tau
\\
&+\int_0^t\int_0^{\alpha_0(t)} g(y,\tau) dy d\tau
-\int_0^{{\alpha_0(t)}} \left( \dfrac{u_0}{\eta_0}\right)(y)dy.
\end{split}
\end{equation}
By virtue of $(\ref{lnv-lnv_0})$, we integrate $(\ref{lnv_xt})$ over $[\alpha_0(t),x]\times[0,t]$ to deduce
\begin{align*}
\ln \dfrac{v(x,t)}{v_0(x)}
= &\ln v(\alpha_0(t),t)-\ln v_0(\alpha_0(t)) + \int_{\alpha_0(t)}^x\left[\left( \dfrac{u}{\eta} \right)(y,t)-\left( \dfrac{u_0}{\eta_0} \right)(y)\right]dy \nonumber
\\[0.5em]
& +\int_0^t\int_{\alpha_0(t)}^x g(y,\tau) dy d\tau + \int_0^t \left[\left( \dfrac{\theta}{\eta v}\right)(x,\tau) - \left( \dfrac{\theta}{\eta v}\right)(\alpha_0(t),\tau) \right]d\tau \nonumber
\\[0.5em]
& + \int_0^t \left[ \left(\dfrac{\chi_x^2}{2 \eta v^2}\right)(x,\tau)-\left(\dfrac{\chi_x^2}{2 \eta v^2}\right)(\alpha_0(t),\tau) \right]d\tau \nonumber
\\[0.5em]
=&\int_0^t \left( \dfrac{\theta}{\eta v} \right)(x,t)d\tau
+ \int_0^t \left(\dfrac{\chi_x^2}{2\eta v^2} \right)(x,t)d\tau + \int_0^t\int_0^x g(y,\tau) dyd\tau \nonumber\\[0.5em]
& - \int_0^t\int_0^1 v \left( \int_0^x g(y,\tau) dy \right)dxd\tau -\int_0^t \int_0^1 \left( \dfrac{\theta+u^2}{\eta}+ \dfrac{\chi_x^2}{2\eta v}\right)dx d\tau \nonumber
\\[0.5em]
& + \int_{\alpha_0(t)}^x \left[\left( \dfrac{u}{\eta} \right)(y,t)-\left( \dfrac{u_0}{\eta_0} \right)(y)\right]dy+\int_0^1 v_0\left(\int_0^x \dfrac{u_0}{\eta_0}dy \right)dx,
\end{align*}
which implies
\begin{equation}\label{vdebiaodashi}
v(x,t)=A(x,t)B(t)D(x,t).
\end{equation}
Here $B(x,t)$ and $D(x,t)$ are given in $(\ref{B})$ and $(\ref{D})$. Besides,
\begin{equation*}
\begin{split}
A(x,t):=& \exp\left\{\int_0^t \left( \dfrac{\theta}{\eta v} \right)(x,t)d\tau+ \int_0^t \left(\dfrac{\chi_x^2}{2\eta v^2} \right)(x,t)d\tau + \int_0^t\int_0^x g(y,\tau) dyd\tau \right.\\[0.5em]
&\left. - \int_0^t\int_0^1 v \left( \int_0^x g(y,\tau) dy \right)dxd\tau \right\}.
\end{split}
\end{equation*}
Noting that
\begin{equation}\label{dA}
\dfrac{d}{dt}A(x,t)=A(x,t)J(x,t)=\dfrac{v(x,t)J(x,t)}{B(t)D(x,t)},
\end{equation}
where $J(x,t)$ is given in $(\ref{J})$, integrating $(\ref{dA})$ over $(0,t)$ yields
\begin{equation*}
A(x,t)=1+\int_0^t\dfrac{v(x,\tau)J(x,\tau)}{B(\tau)D(x,\tau)}d\tau.
\end{equation*}
By inserting the above equality into $(\ref{vdebiaodashi})$, we derive $(\ref{Lv})$.
\hfill$\Box$

\begin{Lemma}\label{Lmaxminv}
There exist two positive constants $C_0$ and $\varepsilon_1$, depending only on $\beta$, $V_0$ and $M_0$,
such that if $(v,u,\chi,\theta)\in X(0,T;m_1,m_2,m_3,N)$ is a solution of the problem $(\ref{L-NSAC1D})$--$(\ref{BC})$
on $(0, T)$, satisfying
\begin{align}
\label{xyjs}
{m_2^{-\alpha}\le2},\quad {(2N)^\alpha\le 1}, \quad {\alpha H(m_1, m_2, m_3, N)\le\varepsilon_1},
\end{align}
with $H(m_1, m_2, m_3, N)\triangleq(1+m_1^{-1}+m_2^{-1}+m_3^{-1}+N)^8$, then
\begin{align}
\label{vsxj}
{C_0\le v(x, t)\le C_0^{-1}}, \qquad \forall(x,t)\in[0,1]\times[0,T].
\end{align}
\end{Lemma}
\noindent{\it\bfseries Proof.}\quad
From $(\ref{energy})$, we get
\begin{equation*}
\int_0^1\left( \dfrac{u^{2}}{2}+\dfrac{(\chi^2-1)^2}{4}+\dfrac{\chi_{x}^2}{2v}+\theta -\ln \theta -1\right)dx
\le E_0.
\end{equation*}
Thanks to $(\ref{v=1})$ and Jessen's inequality, we have
\begin{equation*}
-\ln \bar{\theta}= -\ln \left( \int_{0}^{1} \theta \, dx\right)\le -\int_0^1 \ln \theta \,dx\le E_0.
\end{equation*}
It implies that
\begin{equation*}
\bar{\theta}\ge e^{-E_0}=: \gamma_1 \in (0,1).
\end{equation*}
Thus, we arrive at
\begin{equation}\label{theta}
\bar{\theta}:=\int_{0}^{1} \theta\, dx\in [\gamma_1,1].
\end{equation}

First, we estimate $D(x,t)$. From $(\ref{v=1})$ and $(\ref{xyjs})$, it holds that
\begin{equation*}
\left|\int_{\alpha_0(t)}^x \dfrac{u}{\eta}dy\right|\le \int_0^1 \dfrac{|u|}{\chi^\alpha}dy\le 2\,||u||_{L^2}\le C.
\end{equation*}
Thus, we have
\begin{equation}\label{Ddesxj}
C^{-1}\le D(x,t) \le C,\quad \forall (x,t)\in [0,1]\times[0,T].
\end{equation}

Next, we estimate $B(t)$ by using $(\ref{v=1})$ and $(\ref{theta})$. It follows from the Sobolev's inequality that
\begin{equation*}
||(\chi-\bar{\chi})(t)||_{L^{\infty}}\le ||\chi_x(t)||_{L^2}\le N,\quad \forall t \in [0,T].
\end{equation*}
In addition, we have
\begin{equation*}
\bar{\chi}\le \int_0^1 |\chi|\,dx \le \int_0^1 \chi^2\,dx + 1 = \int_0^1 (\chi^2-1)\,dx + 2 \le \int_0^1 (\chi^2-1)^2 \,dx + 3 \le 7,\quad \forall t \in [0,T].
\end{equation*}
Thus, we show that $||\chi||_{L^{\infty}(Q_T)}\le 7+N\le 2N$. Noting that
\begin{align*}
&\int_0^1 \left( \dfrac{\theta+u^2}{\eta}+ \dfrac{\chi_x^2}{2\eta v}\right)dx\le 2 m_2^{-\alpha} \int_0^1 \left( \theta+\dfrac{u^2}{2}+ \dfrac{\chi_x^2}{2v}\right)dx \le 4  ,
\\
&\int_0^1 \left( \dfrac{\theta+u^2}{\eta}+ \dfrac{\chi_x^2}{2\eta v}\right)dx\ge (2N)^{-\alpha} \int_0^1 \left( \theta+\dfrac{u^2}{2}+ \dfrac{\chi_x^2}{2v}\right)dx \ge \int_0^1 \theta \, dx \ge \gamma_1 ,
\end{align*}
we obtain
\begin{equation}\label{B_}
e^{-4 t} \le B(t) \le e^{-\gamma_1 t}.
\end{equation}
Furthermore, we get
\begin{equation*}
e^{-4 (t-\tau)} \le \dfrac{B(t)}{B(\tau)} \le e^{-\gamma_1 (t-\tau)}.
\end{equation*}
In terms of the definition of $g$, by $(\ref{xyjs})$, we get
\begin{equation}
\begin{split}\label{vg}
&\quad\left|v\int_{0}^x g \, dy\right|\le ||v||_{L^{\infty}}\int_0^1 \left| u \left( \dfrac{1}{\eta}\right)_{t}+ \dfrac{\theta}{v}  \left( \dfrac{1}{\eta}\right)_{x} + \dfrac{\chi_{x}^2}{2v^2} \left( \dfrac{1}{\eta}\right)_{x} + \dfrac{\eta_x u_x}{\eta v} \right|dx \\
& \le \alpha  ||v||_{L^{\infty}} \int_0^1 \left( |\chi^{-\alpha-1}\chi_t u| + \left|\chi^{-\alpha-1} \dfrac{\theta}{v}\chi_x \right| + \left|\chi^{-\alpha-1} \dfrac{\chi_x^2}{2v^2}\chi_x \right| + \left|\dfrac{u_x\chi_x}{v\chi}\right| \right)dx \\
& \le 2 \alpha N m_2^{-\alpha}\left( \dfrac{1}{m_2}||\chi_t||_{L^2}||u||_{L^2} + \dfrac{1}{m_1 m_2}||\theta||_{L^2}||\chi_x||_{L^2}+ \dfrac{1}{m_1^2 m_2}||\chi_x||_{L^\infty}^3 \right) +  \dfrac{2\alpha N}{m_1m_2}||u_x||_{L^2}||\chi_x||_{L^2} \\ 
& \le C\alpha H(m_1,m_2,m_3,N)+ \alpha N ||\chi_t||_{L^2}^2,
\end{split}
\end{equation}
where we have used Cauchy-Schwartz's inequality and the following facts
\begin{align*}
&||v||_{L^{\infty}} \le ||v||_{L^1} + ||v_x||_{L^2}\le 1+ N\le 2N,\quad ||\chi_x||_{L^\infty}^3 \le ||\chi_{xx}||_{L^2}^3 \le N^3,
\\
&||\theta||_{L^2} \le ||\theta-\bar{\theta}||_{L^2} + ||\bar{\theta}||_{L^2}\le||\theta_x||_{L^2}+ 1 \le  N + 1 \le 2N,\quad ||u||_{L^2}\le||u_x||_{L^2}\le N.
\end{align*}
In a similar manner, we can deduce that
\begin{equation}\label{vvg}
\left|v\int_0^1 v \left(\int_{0}^x g \, dy \right) dx\right|\le ||v||_{L^{\infty}}^2 \int_{0}^1 |g| \, dx \le C\alpha H(m_1,m_2,m_3,N)+ \alpha N ||\chi_t||_{L^2}^2.
\end{equation}
Let $f_{+}:=\max\{f,0\}$. Thanks to the mean value theorem, $(\ref{v=1})$ and $(\ref{theta})$, we have
\begin{equation*}
\left( \bar{\theta}^{\frac{\beta+1}{2}}(t) - \theta^{\frac{\beta+1}{2}}(x,t) \right)_+ \le C \left( \int_0^1 \dfrac{\theta^\beta \theta_x^2}{v\theta^2}dx \right)^{\frac{1}{2}}\left( \int_0^1 v\theta \chi_{(\theta\le \bar{\theta})}dx \right)^{\frac{1}{2}}\le CW^{\frac{1}{2}}(t),
\end{equation*}
where $\chi_{(\theta\le \bar{\theta})}=1$ if $\theta\le \bar{\theta}$, and $\chi_{(\theta\le \bar{\theta})}=0$ if $\theta > \bar{\theta}$.
Then by Young's inequality, we get
\begin{equation}\label{min theta}
\min\limits_{x\in[0,1]}\theta(x,t)\ge C_1 - C_2 W(t).
\end{equation}

Using $(\ref{energy})$, $(\ref{Ddesxj})$, $(\ref{B_})$--$(\ref{vg})$, $(\ref{vvg})$ and $(\ref{min theta})$, we infer from $(\ref{Lv})$ and $(\ref{xyjs})$ that
\begin{equation*}
\begin{split}
v(x,t) & \ge  C^{-1} \int_0^t e^{-4(t-\tau)}\min\limits_{x\in[0,1]}\theta \,d\tau - C\alpha N \int_0^t e^{-\gamma_1(t-\tau)}||\chi_t||_{L^2}^2 \,d\tau \\
& - \int_0^t C\alpha e^{-\gamma_1(t-\tau)}H(m_1,m_2,m_3,N)\,d\tau \\
& \ge C^{-1}\int_0^t e^{-4(t-\tau)}[C_1 - C_2 W(\tau)]\,d\tau - C\alpha H(m_1,m_2,m_3,N) \\
& \ge \dfrac{C_1}{4C}(1-e^{-4 t})-\dfrac{C_2}{C}\int_0^t e^{-4(t-\tau)} W(\tau)\,d\tau - C \varepsilon_1.
\end{split}
\end{equation*}
In view of $(\ref{energy})$, we get
\begin{equation*}
\begin{split}
\int_0^t e^{-4(t-\tau)} W(\tau)\,d\tau & = \int_0^{\frac{t}{2}} e^{-4(t-\tau)} W(\tau)\,d\tau +  \int_{\frac{t}{2}}^t e^{-4(t-\tau)} W(\tau)\,d\tau \\
& \le e^{-2t} \int_0^{\frac{t}{2}} W(\tau)\,d\tau + \int_{\frac{t}{2}}^t  W(\tau)\,d\tau \\
& \to 0, \quad {\rm as} \quad t\to \infty.
\end{split}
\end{equation*}
Hence, we can choose a $\tilde{T}$ sufficiently large such that
\begin{equation*}
 \dfrac{C_1e^{-4 t}}{4C}+ \dfrac{C_2}{C}\int_0^t e^{-4(t-\tau)} W(\tau)\,d\tau \le \dfrac{C_1}{16C},\quad t\ge \tilde{T}.
\end{equation*}
Then we get
\begin{equation}\label{vdexj_1}
 v(x,t) \ge \dfrac{C_1}{8C},\quad \forall x \in [0,1], \, t\ge \tilde{T},
\end{equation}
provided $\varepsilon_1>0$ is chosen to be small enough such that $\varepsilon_1\le \min{\left\{ 1,{C_1}/{(16C)} \right\}}$.

For $(x,t)\in [0,1]\times [0,\tilde{T}]$, by $(\ref{xyjs})$, $(\ref{Ddesxj})$, $(\ref{B_})$--$(\ref{vg})$ and $(\ref{vvg})$, we derive from  $(\ref{Lv})$ that
\begin{equation}\label{vdexj_2}
\begin{split}
v(x,t) & \ge B(t)D(x,t) - C\alpha \int_0^{\tilde{T}} e^{-\gamma_1(\tilde{T}-\tau)}[N||\chi_t||_{L^2}^2 + H(m_1,m_2,m_3,N)]\,d\tau \\
& \ge C_3^{-1}e^{-4\tilde{T}}- C_3\alpha H(m_1,m_2,m_3,N) \\
& \ge C_3^{-1}e^{-4\tilde{T}}- C_3 \varepsilon_1 \ge \dfrac{e^{-4\tilde{T}}}{2C_3},
\end{split}
\end{equation}
provided $\varepsilon_1$ is chosen to be such that $\varepsilon_1 \le e^{-4\tilde{T}}/{(2C_3^2)}$.
Combining $(\ref{vdexj_1})$ and $(\ref{vdexj_2})$ gives
\begin{equation}\label{vdexj}
 v(x,t) \ge C_0 := \min \left\{ \dfrac{C_1}{8C}, \dfrac{e^{-4\tilde{T}}}{2C_3} \right\},
\end{equation}
provided $\alpha H(m_1,m_2,m_3,N)\le \varepsilon_1$ with $\varepsilon_1\le \min\{1,{C_1}/{(16C)},e^{-4\tilde{T}}/{(2C_3^2)}\}$.

In what follows, we will deduce the upper bounds of $v$.
Combining $v$ $(\ref{Lv})$ with $(\ref{Ddesxj})$, $(\ref{B_})$--$(\ref{vvg})$, we have
\begin{equation}\label{maxvzhongjian}
\begin{split}
v(x,t) & \le C + C \int_0^t e^{-\gamma_1(t-\tau)}\left( \max\limits_{x\in [0,1]}\theta(x,\tau)+ \max\limits_{x\in [0,1]}\left( \dfrac{\chi_x}{v}\right)^2 \max\limits_{x\in [0,1]}v(x,\tau) \right) d\tau \\
& + C\alpha H(m_1,m_2,m_3,N).
\end{split}
\end{equation}
Noting that for $\beta>0$, there holds
\begin{equation*}
\begin{split}
 \theta^{\frac{\beta+1}{2}}(x,t) - \bar{\theta}^{\frac{\beta+1}{2}}(t) \le  C \left( \int_0^1 \dfrac{\theta^\beta \theta_x^2}{v\theta^2}dx \right)^{\frac{1}{2}}\left( \int_0^1 v\theta \, dx \right)^{\frac{1}{2}}\le CW^{\frac{1}{2}}(t) \max\limits_{x\in [0,1]}v^{\frac{1}{2}}(x,t).
\end{split}
\end{equation*}
Thus, by Young's inequality, we get
\begin{equation}\label{maxtheta}
\theta(x,t) \le C + C W(t)\max\limits_{x\in [0,1]}v(x,t).
\end{equation}
Moreover, by $(\ref{energy})$, for any $(x,t)\in [0,1]\times[0,T]$, it holds that
\begin{equation}\label{chidesj}
\begin{split}
|\chi(x,t)| &\le \left| \int_0^1 (\chi(x,t)-\chi(y,t))dy \right| + \left| \int_0^1 \chi(y,t) dy \right|  \\
& \le \left| \int_0^1 \left(\int_y^x\chi_x(\xi,t)d\xi\right) dy\right| + C \\
& \le \int_0^1 |\chi_x|\,dy + C \le \left(\int_0^1 \dfrac{\chi_x^2}{v}\, dy\right)^{\frac{1}{2}}\left( \int_0^1 v \, dy \right)^{\frac{1}{2}} + C \\
& \le C.
\end{split}
\end{equation}
In view of  $(\ref{L-NSAC1D})_4$, $(\ref{energy})$, $(\ref{vdexj})$ and $(\ref{chidesj})$, we derive that
\begin{equation*}
\begin{split}
\max\limits_{x\in [0,1]}\left( \dfrac{\chi_x}{v}\right)^2 &\le C \int_0^1 \dfrac{\chi_x}{v}\left(\dfrac{\chi_x}{v}\right)_x dx \\
& = C \int_0^1 \dfrac{\chi_x}{v}\left[\left(\dfrac{\chi_x}{v}\right)_x - (\chi^3-\chi)+(\chi^3-\chi)\right]dx \\
& \le C \int_0^1 \max\limits_{x\in [0,1]}\left( \dfrac{\chi_x}{v}\right) |\mu| \,dx + C \int_0^1 \left| \left( \dfrac{\chi_x}{v}\right)(\chi^3-\chi) \right| dx \\
& \le \tilde{\varepsilon} \max\limits_{x\in [0,1]}\left( \dfrac{\chi_x}{v}\right)^2 \int_0^1 \dfrac{\theta}{v}\,dx + C\int_0^1 \dfrac{v\mu^2}{\theta}dx + C \\
& \le \tilde{\varepsilon}C_0^{-1} \max\limits_{x\in [0,1]}\left( \dfrac{\chi_x}{v}\right)^2 + CW(t) + C,
\end{split}
\end{equation*}
where we have used Cauchy-Schwartz's inequality.
We can choose $\tilde{\varepsilon}$ small enough such that
$\tilde{\varepsilon}C_0^{-1}\le 1/2 $, then
\begin{equation}\label{maxchiv}
\max\limits_{x\in [0,1]}\left( \dfrac{\chi_x}{v}\right)^2 \le C + CW(t).
\end{equation}
Combining $(\ref{maxtheta})$ with $(\ref{maxchiv})$, we get
\begin{equation*}
\begin{split}
\max\limits_{x\in [0,1]}\theta(x,\tau)+ & \max\limits_{x\in [0,1]}\left( \dfrac{\chi_x}{v}\right)^2 \max\limits_{x\in [0,1]}v(x,\tau)\\
& \le C + C W(t)\max\limits_{x\in [0,1]}v(x,t)+ (C + CW(t))\max\limits_{x\in [0,1]}v(x,t) \\
& \le C + C(1+W(t))\max\limits_{x\in [0,1]}v(x,t).
\end{split}
\end{equation*}
It follows from $(\ref{maxvzhongjian})$ that
\begin{equation*}
\begin{split}
v(x,t) &\le C + C \int_0^t e^{-\gamma_1(t-\tau)}\left [ C + C(1+W(\tau))\max\limits_{x\in [0,1]}v(x,t) \right] d\tau \\
& \le C + C \int_0^t \left ( e^{-\gamma_1(t-\tau)}+ W(\tau)\right)\max\limits_{x\in [0,1]}v(x,t)\,d\tau.
\end{split}
\end{equation*}
Combining the above inequality with $(\ref{energy})$ and Gronwall's inequality leads to the desired upper bound of specific volume. This, together with $(\ref{vdexj})$, finishes the proof of Lemma $\ref{Lmaxminv}$.
\hfill$\Box$

\vskip2mm
Base on the upper and lower bounds of $v$, we can derive the bounds of $\chi$, one can find the proof in Appendix.
The following lower order a priori estimates are crucial for the high order estimates which require time-independent bounds.

\begin{Lemma}\label{L-chit}
Let the conditions of Lemma $\ref{Lmaxminv}$ be in force. Then there hold
\begin{equation}\label{manyL^1}
||(u_x,\chi_t,\chi^2-1)||_{L^1} + \int_0^1 \dfrac{\chi_x^2}{v}\,dx + \int_0^1 (\chi^2-1)^2\,dx\le CW^{\frac{1}{2}}(t),
\quad t\ge0,
\end{equation}
and
\begin{equation}\label{theta-1 L^infty}
\int_0^T \left|\left|\theta^{\frac{1}{2}}-1\right|\right|_{L^\infty}^2 dt \le C.
\end{equation}
\end{Lemma}

\noindent{\it\bfseries Proof.}\quad
First, in view of  $(\ref{energy})$, $(\ref{v=1})$ and $(\ref{L-NSAC1D})_4$, we have
\begin{equation*}
\begin{split}
||u_x||_{L^1}& = \int_0^1|u_x|\,dx = \int_0^1 \dfrac{\sqrt{\eta(\chi)}\,|u_x|}{\sqrt{v}\sqrt{\theta}}\cdot \dfrac{\sqrt{v}\sqrt{\theta}}{\sqrt{\eta(\chi)}}\,dx \\
& \le C\left( \int_0^1  \dfrac{\eta(\chi)u_x^2}{v\theta} dx\right)^{\frac{1}{2}}\left( \int_0^1 \theta \, dx \right)^{\frac{1}{2}} \\
& \le CW^{\frac{1}{2}}(t),
\end{split}
\end{equation*}
and
\begin{equation}
\begin{split}\label{chi_t L1}
||\chi_t||_{L^1}   = \int_0^1|\chi_t|\,dx & = \int_0^1 \dfrac{|\chi_t|}{\sqrt{v}\sqrt{\theta}}\cdot \sqrt{v}\sqrt{\theta}\, dx \le C\left( \int_0^1  \dfrac{\chi_t^2}{v\theta} dx\right)^{\frac{1}{2}}\left( \int_0^1 \theta \, dx \right)^{\frac{1}{2}}\\
& \le CW^{\frac{1}{2}}(t).
\end{split}
\end{equation}

Next, it follows from $(\ref{L-NSAC1D})_3$ and $(\ref{L-NSAC1D})_4$ that
\begin{equation}\label{1-chi^2 bds}
(1-\chi^2) + \left(\dfrac{\chi_x}{v}\right)_x \dfrac{1}{\chi}= \dfrac{\chi_t}{\chi v}.
\end{equation}
Integrating the above equality over $[0,1]$, by using  $(\ref{chi_t L1})$, we obtain
\begin{equation}\label{1-chi^2}
\int_0^1(1-\chi^2)dx + \int_0^1 \dfrac{\chi_x^2}{v\chi^2}\, dx = \int_0^1 \dfrac{\chi_t}{\chi v}dx\le C||\chi_t||_{L^1}\le CW^{\frac{1}{2}}(t).
\end{equation}
In view of $(\ref{L-NSAC1D})_3$ and $(\ref{L-NSAC1D})_4$, we have
\begin{equation}\label{chi_t/v}
\dfrac{\chi_t}{v} = \left(\dfrac{\chi_x}{v}\right)_x -(\chi^3-\chi).
\end{equation}
Multiplying $(\ref{chi_t/v})$ by $\chi$, and using $(\ref{chi_t L1})$, $(\ref{1-chi^2})$, we get
\begin{equation*}
\begin{split}
\int_0^1 \dfrac{\chi_x^2}{v}\,dx + \int_0^1 (\chi^2-1)^2\,dx & = \int_0^1 (1-\chi^2)\,dx -\int_0^1\dfrac{\chi\chi_t }{v}\,dx \\
& \le  CW^{\frac{1}{2}}(t) + C ||\chi_t||_{L^1} \\
& \le   CW^{\frac{1}{2}}(t).
\end{split}
\end{equation*}
From $(\ref{1-chi^2 bds})$, in a similar manner to $(\ref{1-chi^2})$, there holds
\begin{equation*}
\int_0^1(\chi^2-1)dx = \int_0^1 \dfrac{\chi_x^2}{v\chi^2}\, dx - \int_0^1 \dfrac{\chi_t}{\chi v}dx\le C\int_0^1 \dfrac{\chi_x^2}{v}\,dx +C||\chi_t||_{L^1}\le CW^{\frac{1}{2}}(t).
\end{equation*}
Which together with $(\ref{1-chi^2})$ implies $ ||\chi^2-1||_{L^1} \le CW^{\frac{1}{2}}(t)$.

Finally,
\begin{equation*}
\left|\left|\theta^{\frac{1}{2}}-1\right|\right|_{L^\infty} \le \left|\left|\theta^{\frac{1}{2}}-\bar{\theta}^{\frac{1}{2}}\right|\right|_{L^\infty} + \left|\left|\bar{\theta}^{\frac{1}{2}}-1\right|\right|_{L^\infty}.
\end{equation*}
For the first item on the right side, using $(\ref{v=1})$ and $(\ref{vsxj})$, for any $\beta>0$, there holds
\begin{equation}\label{theta1/2_1}
|\theta^{\frac{1}{2}}-\bar{\theta}^{\frac{1}{2}}| \le \left| \theta^{\frac{\beta+1}{2}}-\bar{\theta}^{\frac{\beta+1}{2}} \right| \le  C \left( \int_0^1 \dfrac{\theta^\beta \theta_x^2}{v\theta^2}dx \right)^{\frac{1}{2}}\left( \int_0^1 v\theta \, dx \right)^{\frac{1}{2}}\le CW^{\frac{1}{2}}(t).
\end{equation}
For the second item, it follows from $(\ref{v=1})$ and $(\ref{manyL^1})$ that
\begin{equation}\label{theta1/2_2}
\begin{split}
|1-\bar{\theta}^{\frac{1}{2}}| & \le  C|1-\bar{\theta}| = C \left| 1 - \int_0^1 \theta \,dx\right| \\
& = C \int_0^1 \dfrac{u^2}{2} dx +  C \int_0^1 \left( \dfrac{(\chi^2-1)^2}{4}+\dfrac{\chi_x^2}{2v} \right) dx\\
& \le C ||u||_{L^\infty}||u||_{L^2} + CW^{\frac{1}{2}}(t) \\
& \le C \int_0^1 |u_x|\,dx +  CW^{\frac{1}{2}}(t) \le  CW^{\frac{1}{2}}(t).
\end{split}
\end{equation}
Combining $(\ref{theta1/2_1})$, $(\ref{theta1/2_2})$ with $(\ref{energy})$ gives
\begin{equation*}
\int_0^T ||\theta^{\frac{1}{2}}-1||_{L^\infty}^2 \,dt  \le  C \int_0^T W(t)\,dt \le C.
\end{equation*}
\hfill$\Box$

\begin{Lemma}\label{Ltheta_x}
Let the conditions of Lemma $\ref{Lmaxminv}$ be in force.
Then for any $p>0$, there exists a positive constant C, which may depend on $p$, such that
\begin{equation}\label{theta^beta theta_x}
\int_0^T \int_0^1 \dfrac{\theta^\beta \theta_x^2}{\theta^{p+1}}\,dxdt \le C(p).
\end{equation}
\end{Lemma}

\noindent{\it\bfseries Proof.}\quad
We only consider the case that $p \ne 1$, since $(\ref{theta^beta theta_x})$ with $p=1$ is thanks to $(\ref{energy})$ and $(\ref{vsxj})$. In fact, multiplying $(\ref{L-NSAC1D})_5$ by $\theta^{-p}$ with $p\ne1$ and integrating by parts, we derive that
\begin{align}\label{Gronwall's_1}
&\quad\dfrac{1}{p-1} \dfrac{d}{dt} \int_0^1 \theta^{1-p}\,dx + p \int_0^1 \dfrac{\theta^\beta \theta_x^2}{ v\theta^{p+1}} \,dx + \int_0^1 \dfrac{\eta(\chi) u_x^2}{ v\theta^{p}} \,dx + \int_0^1 \dfrac{v \mu^2}{\theta^{p}} \,dx  \nonumber \\
& =  \int_0^1 \dfrac{(\theta^{1-p}-1)u_x}{v} \,dx  + \left( \int_0^1 \ln v \,dx \right)_t   \nonumber \\
& \le C(p)\int_0^1 |\theta^{\frac{1}{2}}-1|\left(\theta^{\frac{1}{2}-p}+1 \right)|u_x| \,dx  + \left( \int_0^1 \ln v \,dx \right)_t   \nonumber \\
& \le C(p)  \left|\left| \theta^{\frac{1}{2}}-1 \right|\right|_{L^\infty} \left( \int_0^1 \dfrac{v\theta^{1-p}}{\eta(\chi)}dx \right)^{\frac{1}{2}}\left( \int_0^1  \dfrac{\eta(\chi)u_x^2}{v\theta^p}dx \right)^{\frac{1}{2}}\nonumber \\
& \quad +  C(p)  \left|\left| \theta^{\frac{1}{2}}-1 \right|\right|_{L^\infty}\int_0^1|u_x| \,dx + \left( \int_0^1 \ln v \,dx \right)_t \nonumber \\
& \le \dfrac{1}{2} \left( \int_0^1  \dfrac{\eta(\chi)u_x^2}{v\theta^p}dx \right) +  C(p)  \left|\left| \theta^{\frac{1}{2}}-1 \right|\right|_{L^\infty}^2 \left(1 + \int_0^1 \theta^{1-p}dx  \right) \nonumber \\
& \quad + C(p) \left( \int_0^1 |u_x| \,dx \right)^2 + \left( \int_0^1 \ln v \,dx \right)_t.
\end{align}

{\rm\bfseries Case I.} $p>1$.

By virtue of Gronwall's inequality, $(\ref{manyL^1})$, $(\ref{theta-1 L^infty})$ and the fact that $||\ln v(t)||_{L^1}\le C$ due to $(\ref{vsxj})$, we can infer from $(\ref{Gronwall's_1})$ that $(\ref{theta^beta theta_x})$ holds.

{\rm\bfseries Case II.} $0<p<1$.

Noting that $||\theta^{1-p}||_{L^1}\le C + C||\theta||_{L^1}$, integrating $(\ref{Gronwall's_1})$ over $[0,T]$, by $(\ref{manyL^1})$, $(\ref{theta-1 L^infty})$ and $||\ln v(t)||_{L^1}\le C$, we easily arrive at $(\ref{theta^beta theta_x})$ with $0<p<1$.

Combining {\rm\bfseries Case I.} with {\rm\bfseries Case II.} finishes the proof of Lemma $\ref{Ltheta_x}$.
\hfill$\Box$

\begin{Lemma}\label{L u_x mu L^2}
Let the conditions of Lemma $\ref{Lmaxminv}$ be in force. Then there holds
\begin{equation}\label{mu + u_x L^2}
\int_0^T ||u_x||_{L^2}^2 \,dt + \int_0^T ||\mu||_{L^2}^2 \,dt \le C.
\end{equation}
\end{Lemma}

\noindent{\it\bfseries Proof.}\quad
Integrating $(\ref{L-NSAC1D})_5$ over $[0,1]\times[0,T]$, using $(\ref{v=1})$, $(\ref{vsxj})$, $(\ref{theta-1 L^infty})$ and the boundary condition $(\ref{BC})$, we find that
\begin{align}
&\quad\int_0^T \int_0^1 v \mu^2 \, dxdt  + \int_0^T \int_0^1 \dfrac{\eta(\chi)u_x^2}{v} \,dxdt \nonumber  \\
& = \int_0^1 \theta(x,T) \, dx -  \int_0^1 \theta_0(x) \, dx + \int_0^T \int_0^1 \dfrac{\theta u_x}{v} \,dxdt
 \nonumber\\
& \le C + \int_0^T \int_0^1 \dfrac{(\theta-1) u_x}{v} \,dxdt + \int_0^T \left( \int_0^1 \ln v \,dx \right)_t dt \nonumber\\
& \le C + \dfrac{1}{2}\int_0^T \int_0^1 \dfrac{\eta(\chi)u_x^2}{v} \,dxdt + \int_0^T \int_0^1 \dfrac{(\theta^{\frac{1}{2}}-1)^2(\theta^{\frac{1}{2}}+1)^2}{\eta(\chi)v} \,dxdt  \nonumber  \\
& \le C + \dfrac{1}{2}\int_0^T \int_0^1 \dfrac{\eta(\chi)u_x^2}{v} \,dxdt + C\int_0^T \left(||\theta^{\frac{1}{2}}-1||_{L^\infty}^2 \int_0^1 (\theta+1)\,dx\right)dt \nonumber  \\
& \le C + \dfrac{1}{2}\int_0^T \int_0^1 \dfrac{\eta(\chi)u_x^2}{v} \,dxdt + C\int_0^T ||\theta^{\frac{1}{2}}-1||_{L^\infty}^2\,dt \nonumber  \\
& \le C + \dfrac{1}{2}\int_0^T \int_0^1 \dfrac{\eta(\chi)u_x^2}{v} \,dxdt. \nonumber
\end{align}
It implies that Lemma $\ref{L u_x mu L^2}$ holds.
\hfill$\Box$

\begin{Lemma}\label{L-chixx}
Let the conditions of Lemma $\ref{Lmaxminv}$ be in force. Then it holds that
\begin{equation} \label{(chi_x/v)_x L2}
\int_0^T \left|\left| \left(\dfrac{\chi_x}{v}\right)_x\right|\right|_{L^2}^2 \,dt + \int_0^T ||\chi^2-1||_{L^2}^2 \,dt \le C.
\end{equation}
\end{Lemma}

\noindent{\it\bfseries Proof.}\quad
From $(\ref{L-NSAC1D})_4$ and $(\ref{vsxj})$, we derive that
\begin{align*}
\int_0^1 \left| \left( \dfrac{\chi_x}{v} \right)_x \right|^2 dx & \le 2 \int_0^1 \mu^2 \,dx + 2 \int_0^1 \chi^2(\chi^2-1)^2 \, dx \nonumber \\
& \le C ||\mu||_{L^2}^2 + C\int_0^1 (\chi^2-1)^2 \, dx \nonumber \\
& \le C ||\mu||_{L^2}^2 + C ||\chi^2-1||_{L^\infty} \int_0^1 |\chi^2-1| \, dx \nonumber \\
& \le C ||\mu||_{L^2}^2 + \tilde{\varepsilon}_1 ||\chi^2-1||_{L^\infty}^2 + C ||\chi^2-1||_{L^1}^2 \nonumber \\
& \le C ||\mu||_{L^2}^2 + (C+\tilde{\varepsilon}_1) ||\chi^2-1||_{L^1}^2 + \tilde{\varepsilon}_1 ||2\chi\chi_x||_{L^1}^2 \nonumber \\
& \le C ||\mu||_{L^2}^2 + (C+\tilde{\varepsilon}_1) ||\chi^2-1||_{L^1}^2 + \tilde{\varepsilon}_1 C\left|\left|\dfrac{\chi_x}{v}\right|\right|_{L^2}^2  \nonumber \\
& \le C ||\mu||_{L^2}^2 + (C+\tilde{\varepsilon}_1) ||\chi^2-1||_{L^1}^2 + \tilde{\varepsilon}_1 \tilde{C} \left|\left| \left(\dfrac{\chi_x}{v}\right)_x\right|\right|_{L^2}^2,
\end{align*}
where we have used Cauchy-Schwartz's inequality and the fact $W^{1,1}([0,1])\hookrightarrow L^{\infty}([0,1])$.
We can chose $\tilde{\varepsilon}$ is small enough with $\tilde{\varepsilon}_1 \tilde{C}\le 1/2$,
integrating the above inequality over $[0,T]$ gives
\begin{equation}\label{chixx}
\int_0^T \left|\left| \left(\dfrac{\chi_x}{v}\right)_x\right|\right|_{L^2}^2 \,dt \le C \int_0^T ||\mu||_{L^2}^2\,dt + C \int_0^T ||\chi^2-1||_{L^1}^2\,dt \le C,
\end{equation}
where we have uesd $(\ref{manyL^1})$ and $(\ref{mu + u_x L^2})$.

In view of $(\ref{L-NSAC1D})_4$, we get
\begin{align}
\int_0^T ||\chi^2-1||_{L^2}^2 \,dt & \le C \int_0^T \int_0^1 \dfrac{1}{\chi^2}\mu^2 \,dxdt +  C \int_0^T \int_0^1 \dfrac{1}{\chi^2} \left| \left( \dfrac{\chi_x}{v} \right)_x \right|^2 \,dxdt \nonumber \\
& \le C \int_0^T ||\mu||_{L^2}^2\,dt + \int_0^T \left|\left| \left(\dfrac{\chi_x}{v}\right)_x\right|\right|_{L^2}^2\,dt
\le C. \nonumber
\end{align}
From which and $(\ref{chixx})$, Lemma $\ref{L-chixx}$ holds.
\hfill$\Box$

\begin{Lemma}\label{L-sup v_x}
Let the conditions of Lemma $\ref{Lmaxminv}$ be in force. Then there holds
\begin{equation}\label{sup v_x chi_x}
\sup \limits_{0\le t \le T}||(v_x(t),\chi_x(t))||_{L^2}^2 + \int_0^T \int_0^1 \left( v_x^2 + \chi_{xx}^2+ \theta v_x^2 \right)\,dxdt \le C_2.
\end{equation}
\end{Lemma}

\noindent{\it\bfseries Proof.}\quad
Based on $(\ref{L-NSAC1D})_1$ and $(\ref{L-NSAC1D})_2$, we deduce that
\begin{equation}\label{v_x}
\begin{split}
\left( \dfrac{\eta(\chi)v_x}{v} \right)_t & = \eta_t \dfrac{v_x}{v} +\eta \left( \dfrac{u_x}{v} \right)_x = \left( \dfrac{\eta u_x}{v} \right)_x + \dfrac{\eta_t u_x - \eta_x u_x}{v}  \\
& = u_t + \dfrac{\theta_x}{v} - \dfrac{\theta v_x}{v^2} + \dfrac{1}{2}\left( \dfrac{\chi_x^2}{v^2} \right)_x + \dfrac{\alpha \chi^\alpha}{v\chi}(\chi_t v_x - \chi_x u_x).
\end{split}
\end{equation}
Multiplying $(\ref{v_x})$ by $\eta(\chi)v_x/v$ and integrating by parts over $[0,1]$ yield
\begin{align}\label{v_x L^infty L^2}
\dfrac{1}{2}& \dfrac{d}{dt} \left| \left| \dfrac{\eta(\chi)v_x}{v}  \right|\right|_{L^2}^2 + \int_0^1 \dfrac{\eta(\chi)\theta v_x^2}{v^3}\,dx \nonumber\\
& = \int_0^1 \dfrac{\eta(\chi)u_t v_x}{v} \,dx + \int_0^1 \dfrac{\eta(\chi) v_x \theta_x}{v^2} \,dx + \dfrac{1}{2}\int_0^1 \left( \dfrac{\chi_x^2}{v^2} \right)_x \dfrac{\eta(\chi)v_x}{v}\, dx  \nonumber\\
&\quad + \int_0^1 \dfrac{\alpha \eta^2(\chi_t v_x-\chi_xu_x)v_x}{v^2\chi}\,dx  \nonumber\\
& = \dfrac{d}{dt} \int_0^1 \dfrac{\eta(\chi)uv_x}{v}\,dx - \int_0^1 \left( \eta u (\ln v)_{xt}+\eta_t u (\ln v)_x \right)\,dx + \int_0^1 \dfrac{\eta(\chi) v_x \theta_x}{v^2} \,dx \nonumber\\
& \quad + \dfrac{1}{2}\int_0^1 \left( \dfrac{\chi_x^2}{v^2} \right)_x \dfrac{\eta(\chi)v_x}{v}\, dx + \int_0^1 \dfrac{\alpha \eta^2(\chi_t v_x-\chi_xu_x)v_x}{v^2\chi}\,dx \nonumber\\
& = \dfrac{d}{dt} \int_0^1 \dfrac{\eta(\chi)uv_x}{v}\,dx - \int_0^1 \left[ \left( \dfrac{\eta u u_x}{v} \right)_x - (\eta u)_x \dfrac{u_x}{v}+ \eta_t u (\ln v)_x \right] \,dx + \int_0^1 \dfrac{\eta(\chi) v_x \theta_x}{v^2} \,dx \nonumber\\
& \quad + \dfrac{1}{2}\int_0^1 \left( \dfrac{\chi_x^2}{v^2} \right)_x \dfrac{\eta(\chi)v_x}{v}\, dx + \int_0^1 \dfrac{\alpha \eta^2(\chi_t v_x-\chi_xu_x)v_x}{v^2\chi}\,dx.
\end{align}
Due to the non-slip boundary conditions $u\big|_{x=0,1}=0$ and the fact
\begin{equation*}
\int_0^1 \dfrac{\eta(\chi)uv_x}{v}\,dx \le \dfrac{1}{2} \int_0^1 \left( \dfrac{\eta(\chi)v_x}{v} \right)^2\,dx + C \int_0^1 u^2 \,dx,
\end{equation*}
we derive after integrating $(\ref{v_x L^infty L^2})$ over $[0,T]$ that
\begin{equation}\label{I_1--I_4}
\begin{split}
&\quad \dfrac{1}{4}\sup \limits_{0\le t \le T} \left| \left| \dfrac{\eta(\chi)v_x}{v}  \right|\right|_{L^2}^2 + \int_0^T\int_0^1 \dfrac{\eta(\chi)\theta v_x^2}{v^3}\,dxdt \\
& \le C + \int_0^T \int_0^1 \left[ (\eta u)_x \dfrac{u_x}{v}-\eta_t u (\ln v)_x \right] \,dxdt + \int_0^T \int_0^1 \dfrac{\eta(\chi) v_x \theta_x}{v^2} \,dxdt \\
& + \dfrac{1}{2} \int_0^T\int_0^1 \left( \dfrac{\chi_x^2}{v^2} \right)_x \dfrac{\eta(\chi)v_x}{v}\, dxdt +  \int_0^T\int_0^1 \dfrac{\alpha \eta^2(\chi_t v_x-\chi_xu_x)v_x}{v^2\chi}\,dxdt \\
& :=C + \sum \limits_{i=1}^{4}I_i.
\end{split}
\end{equation}
It follows from $(\ref{xyjs})$ and $(\ref{vsxj})$ that
\begin{equation}\label{I_1}
\begin{split}
I_1 & =  \int_0^T \int_0^1 \left[ (\eta u)_x \dfrac{u_x}{v}-\eta_t u (\ln v)_x \right] \,dxdt \\
& \le C \int_0^T \int_0^1 \left[ \eta u_x^2 + \alpha \chi^{\alpha-1}\left( |\chi_x u u_x| + |\chi_t u v_x| \right) \right] \,dxdt \\
& \le C + C\alpha \int_0^T  \left[ ||u_x||_{L^2}^2 \left( ||u_x||_{L^2}^2 + ||v_x||_{L^2}^2 \right) + ||\chi_x||_{L^2}^2 + ||\chi_t||_{L^2}^2\right] \,dt \\
& \le C + C\alpha N^2 \int_0^T  \left( ||u_x||_{L^2}^2 + \left|\left|\dfrac{\chi_x}{v}\right|\right|_{L^2}^2 + ||\chi_t||_{L^2}^2 \right) \,dt \\
& \le C + C\alpha N^2 \int_0^T  \left( ||u_x||_{L^2}^2 + \left|\left|\left(\dfrac{\chi_x}{v}\right)_x\right|\right|_{L^2}^2 + ||\chi_t||_{L^2}^2 \right) \,dt  \\
& \le C + C\alpha H(m_1,m_2,m_3,N) \le C,
\end{split}
\end{equation}
where we have used the fact $||u||_{L^\infty} \le C ||u||_{L^2}+ C||u_x||_{L^2}\le C||u_x||_{L^2}$.
We choose $p=\beta>0$ in $(\ref{theta^beta theta_x})$ and obtain
\begin{align}\label{theta_x/theta}
\int_0^T \int_0^1 \dfrac{\theta_x^2}{\theta}\,dxdt \le C(\beta)\le C.
\end{align}
Using $(\ref{vsxj})$ and $(\ref{theta_x/theta})$ gives
\begin{equation*}
\begin{split}
I_2 & = \int_0^T \int_0^1 \dfrac{\eta(\chi) v_x \theta_x}{v^2} \,dxdt \\
& \le \dfrac{3}{8} \int_0^T\int_0^1 \dfrac{\eta(\chi)\theta v_x^2}{v^3}\,dxdt + C \int_0^T\int_0^1 \dfrac{\eta(\chi)\theta_x^2}{v \theta} \,dxdt \\
& \le \dfrac{3}{8} \int_0^T\int_0^1 \dfrac{\eta(\chi)\theta v_x^2}{v^3}\,dxdt + C \int_0^T\int_0^1 \dfrac{\theta_x^2}{\theta} \,dxdt \\
& \le \dfrac{3}{8} \int_0^T\int_0^1 \dfrac{\eta(\chi)\theta v_x^2}{v^3}\,dxdt + C.
\end{split}
\end{equation*}
In view of $(\ref{(chi_x/v)_x L2})$, we have
\begin{equation*}
\begin{split}
I_3 & = \dfrac{1}{2} \int_0^T\int_0^1 \left( \dfrac{\chi_x^2}{v^2} \right)_x \dfrac{\eta(\chi)v_x}{v}\, dxdt \\
& \le \int_0^T\int_0^1 \left( \dfrac{\chi_x}{v}\right)^2 \left( \dfrac{\eta(\chi)v_x}{v}\right)^2dxdt + \int_0^T\int_0^1 \left|\left( \dfrac{\chi_x}{v}\right)_x \right|^2\,dxdt  \\
& \le \int_0^T \left|\left| \dfrac{\chi_x}{v}\right|\right|_{L^\infty}^2\int_0^1 \left( \dfrac{\eta(\chi)v_x}{v}\right)^2dxdt + C  \\
& \le \int_0^T \left|\left| \left( \dfrac{\chi_x}{v}\right)_x\right|\right|_{L^2}^2  \left|\left| \dfrac{\eta(\chi)v_x}{v}\right|\right|_{L^2}^2dt + C.
\end{split}
\end{equation*}
In terms of $(\ref{xyjs})$, $(\ref{vsxj})$ and $(\ref{(chi_x/v)_x L2})$, then we get
\begin{equation}\label{I_4}
\begin{split}
I_4 & = \int_0^T\int_0^1 \dfrac{\alpha \eta^2(\chi_t v_x-\chi_xu_x)v_x}{v^2\chi}\,dxdt \\
& \le C \alpha \int_0^T \int_0^1 \sqrt{\eta}\left(v_x^2|\chi_t|+ |\chi_x u_x v_x | \right)dxdt \\
& \le \dfrac{3}{8} \int_0^T \int_0^1 \dfrac{\eta(\chi)\theta v_x^2}{v^3}\,dx dt  + C\alpha^2 \int_0^T \int_0^1 \dfrac{v^3}{\theta}\left( v_x^2 \chi_t^2 + \chi_x^2 u_x^2 \right)dx dt \\
& \le \dfrac{3}{8} \int_0^T \int_0^1 \dfrac{\eta(\chi)\theta v_x^2}{v^3}\,dx dt  + C\dfrac{\alpha^2 N^2}{m_3}\int_0^T \int_0^1 \left( \chi_t^2 + u_x^2 \right)dx dt \\
& \le \dfrac{3}{8} \int_0^T \int_0^1 \dfrac{\eta(\chi)\theta v_x^2}{v^3}\,dx dt  + C[\alpha H(m_1,m_2,m_3,N)]^2 \\
& \le \dfrac{3}{8} \int_0^T \int_0^1 \dfrac{\eta(\chi)\theta v_x^2}{v^3}\,dx dt  + C,
\end{split}
\end{equation}
where we have uesd the fact $||(v_x,\chi_x)||_{L^\infty}\le C||(v_x,\chi_x)||_{H^1}\le CN$ for $t\in[0,T]$.

Substituting $(\ref{I_1})$--$(\ref{I_4})$ into $(\ref{I_1--I_4})$, by virtue of Gronwall's inequality,
$(\ref{(chi_x/v)_x L2})$ and $(\ref{vsxj})$, we deduce that
\begin{equation*}
\sup \limits_{0\le t \le T}||v_x(t)||_{L^2}^2 + \int_0^T \int_0^1 \theta v_x^2 \,dxdt \le C.
\end{equation*}
From which and $(\ref{theta-1 L^infty})$, we get
\begin{equation*}
\int_0^T \int_0^1v_x^2 \,dxdt \le C \int_0^T  \left|\left|\theta^{\frac{1}{2}}-1\right|\right|_{L^\infty}^2 ||v_x||_{L^2}^2 \,dt + C \int_0^T \int_0^1 \theta v_x^2 \,dxdt \le C.
\end{equation*}

We rewrite $(\ref{L-NSAC1D})_{3,4}$ as follows
\begin{equation}\label{chi_t=chi_xx}
\chi_t = v \left( \dfrac{\chi_x}{v} \right)_x - v(\chi^3-\chi)=\chi_{xx}-\dfrac{\chi_x v_x}{v}-v(\chi^3-\chi).
\end{equation}
Multiplying $(\ref{chi_t=chi_xx})$ by $\chi_{xx}$ and integrating the result over $[0,1]$ yield
\begin{equation*}
\begin{split}
\dfrac{1}{2}\dfrac{d}{dt}||\chi_x||_{L^2}^{2} + ||\chi_{xx}||_{L^2}^{2} & = \int_0^1 \dfrac{\chi_x v_x}{v}\chi_{xx}\,dx + \int_0^1 v(\chi^3-\chi)\chi_{xx}\,dx \\
& \le \dfrac{1}{2}||\chi_{xx}||_{L^2}^{2} + C \int_0^1 \chi_x^2 v_x^2\,dx + C \int_0^1 \chi^2(\chi^2-1)^2\,dx\\
& \le \dfrac{1}{2}||\chi_{xx}||_{L^2}^{2} +  C \left|\left|\dfrac{\chi_x}{v}\right|\right|_{L^\infty}^2 ||v_x||_{L^2}^{2} + C ||\chi^2-1||_{L^2}^{2} \\
& \le \dfrac{1}{2}||\chi_{xx}||_{L^2}^{2} +  C \left|\left|\left(\dfrac{\chi_x}{v}\right)_x\right|\right|_{L^2}^2  + C ||\chi^2-1||_{L^2}^{2}.
\end{split}
\end{equation*}
By virtue of Gronwall's inequality and $(\ref{(chi_x/v)_x L2})$, we arrive at
\begin{equation*}
\sup \limits_{0\le t \le T}||\chi_x||_{L^2}^2 + \int_0^T \int_0^1 \chi_{xx}^2\,dxdt \le C.
\end{equation*}
This completes the proof of Lemma $\ref{L-sup v_x}$.
\hfill$\Box$

\begin{Lemma}
Let the conditions of Lemma $\ref{Lmaxminv}$ be in force. Then we have
\begin{equation}\label{sup chi_xx chi_t}
\sup \limits_{0\le t \le T}\left|\left| \left( \chi_{xx},\chi_t,\left(\dfrac{\chi_x}{v}\right)_x\,\right) \right|\right|_{L^2}^2 + \int_0^T \left|\left| \left( \chi_{xt},\left(\dfrac{\chi_x}{v}\right)_t\,\right) \right|\right|_{L^2}^2 dt \le C_3.
\end{equation}
\end{Lemma}

\noindent{\it\bfseries Proof.}\quad
Rewrite $(\ref{L-NSAC1D})_{3,4}$ as
\begin{equation}\label{chi_t/v}
\dfrac{\chi_t}{v} =  \left( \dfrac{\chi_x}{v} \right)_x - (\chi^3-\chi).
\end{equation}
Differentiating $(\ref{chi_t/v})$ with respect to $x$ yields
\begin{equation*}
\left(\dfrac{\chi_x}{v}\right)_t - \left( \dfrac{\chi_x}{v} \right)_{xx} = \dfrac{\chi_t v_x}{v^2 }- \dfrac{\chi_x u_x}{v^2} - \chi_x(3\chi^2-1).
\end{equation*}
Multiplying the above equation by $\left(\frac{\chi_x}{v}\right)_t $ and integrating over $[0,1]$, we get
\begin{equation*}
\begin{split}
&\dfrac{1}{2} \dfrac{d}{dt}\left| \left|\left(\dfrac{\chi_x}{v}\right)_x \right|\right|_{L^2}^{2}  + \left|\left|\left(\dfrac{\chi_x}{v}\right)_t \right|\right|_{L^2}^{2}
\\
= &\int_0^1\left( \dfrac{\chi_t v_x}{v^2 }- \dfrac{\chi_x u_x}{v^2} - \chi_x(3\chi^2-1) \right)\left(\dfrac{\chi_x}{v}\right)_t \,dx
\\
\le& \dfrac{1}{4}\left|\left|\left(\dfrac{\chi_x}{v}\right)_t \right|\right|_{L^2}^{2} + C \int_0^1 \left( \chi_t^2 v_x^2 + \chi_x^2 u_x^2 + \chi_x^2(3\chi^2-1)^2 \right)\,dx
\\
\le& \dfrac{1}{4}\left|\left|\left(\dfrac{\chi_x}{v}\right)_t \right|\right|_{L^2}^{2} + C ||\chi_{t}||_{L^\infty}^{2}||v_{x}||_{L^2}^{2} + C ||\chi_{x}||_{L^\infty}^{2}||u_{x}||_{L^2}^{2} + C ||\chi_{x}||_{L^2}^{2}
\\
\le& \dfrac{1}{4}\left|\left|\left(\dfrac{\chi_x}{v}\right)_t \right|\right|_{L^2}^{2} + C ||\chi_t||_{L^2}||\chi_{xt}||_{L^2} + C ||\chi_{xx}||_{L^2}^{2}
\\
= &\dfrac{1}{4}\left|\left|\left(\dfrac{\chi_x}{v}\right)_t \right|\right|_{L^2}^{2} + C ||\chi_t||_{L^2}\left|\left|\left[\left(\dfrac{\chi_x}{v}\right)_t + \dfrac{\chi_x u_x}{v^2} \right] v \right|\right|_{L^2} + C ||\chi_{xx}||_{L^2}^{2}
\\
\le &\dfrac{1}{2}\left|\left|\left(\dfrac{\chi_x}{v}\right)_t \right|\right|_{L^2}^{2} + C ||\chi_t||_{L^2}^2 + C ||\chi_{x}||_{L^\infty}^{2}||u_{x}||_{L^2}^{2} + C ||\chi_{xx}||_{L^2}^{2}
\\
\le& \dfrac{1}{2}\left|\left|\left(\dfrac{\chi_x}{v}\right)_t \right|\right|_{L^2}^{2} + C ||\chi_t||_{L^2}^2 + C ||\chi_{xx}||_{L^2}^{2}.
\end{split}
\end{equation*}
Integrating it over $[0,T]$ gives
\begin{equation*}
\sup \limits_{0\le t \le T}\left|\left| \left(\dfrac{\chi_x}{v}\right)_x\right|\right|_{L^2}^2 + \int_0^T \left|\left|\left(\dfrac{\chi_x}{v}\right)_t \right|\right|_{L^2}^2 dt \le C.
\end{equation*}
Thus, it follows from $(\ref{vsxj})$ and the above inequality that
\begin{equation*}
\begin{split}
||\chi_{xx}||_{L^2}^{2}= \left|\left|\left[\left(\dfrac{\chi_x}{v}\right)_x + \dfrac{\chi_x v_x}{v^2} \right] v \right|\right|_{L^2}^{2} & \le C \left|\left|\left(\dfrac{\chi_x}{v}\right)_x \right|\right|_{L^2}^{2} + C \left|\left|\dfrac{\chi_x}{v}\right|\right|_{L^\infty}^{2}||v_x||_{L^2}^2 \\
& \le C \left|\left|\left(\dfrac{\chi_x}{v}\right)_x \right|\right|_{L^2}^{2} \le C.
\end{split}
\end{equation*}
Moreover,
\begin{equation*}
\begin{split}
||\chi_{t}||_{L^2}^{2}= \left|\left| v\left(\dfrac{\chi_{x}}{v}\right)_x-v(\chi^3-\chi)\right|\right|_{L^2}^{2}\le  C \left|\left|\left(\dfrac{\chi_x}{v}\right)_x \right|\right|_{L^2}^{2} + C \le C.
\end{split}
\end{equation*}
Similarly, we can deduce that
\begin{equation*}
\begin{split}
\int_0^T||\chi_{xt}||_{L^2}^{2}\,dt & = \int_0^T \left|\left|\left[\left(\dfrac{\chi_x}{v}\right)_t + \dfrac{\chi_x u_x}{v^2} \right] v \right|\right|_{L^2}^{2}dt \\
& \le C \int_0^T \left|\left|\left(\dfrac{\chi_x}{v}\right)_t \right|\right|_{L^2}^{2}\,dt + C \int_0^T\left|\left|\dfrac{\chi_x}{v}\right|\right|_{L^\infty}^{2}||u_x||_{L^2}^2 \,dt \le C.
\end{split}
\end{equation*}
Then we see that $(\ref{sup chi_xx chi_t})$ holds.
\hfill$\Box$

\begin{Lemma}
Let the conditions of Lemma $\ref{Lmaxminv}$ be in force, then
\begin{equation}\label{sup u_x 00}
\sup \limits_{0\le t \le T}\left|\left| u_x \right|\right|_{L^2}^2 + \int_0^T \left|\left| \left(u_t,u_{xx}, \theta_x \right) \right|\right|_{L^2}^2 dt \le C_4.
\end{equation}
\end{Lemma}

\noindent{\it\bfseries Proof.}\quad
We rewrite $(\ref{L-NSAC1D})_{2}$ as follows
\begin{equation*}
u_t - \dfrac{\eta(\chi)u_{xx}}{v} = \dfrac{\eta_x u_{x}}{v} - \dfrac{\eta u_{x}v_x}{v^2} - \dfrac{\theta_{x}}{v} + \dfrac{\theta v_x}{v} - \dfrac{1}{2}\left( \dfrac{\chi_x^2}{v^2} \right)_x.
\end{equation*}
Multiplying it by $u_{xx}$, then integrating the result over $[0,1]\times[0,T]$, by using $(\ref{vsxj})$ we have
\begin{align}\label{u_x 0}
&\quad \dfrac{1}{2} \dfrac{d}{dt} ||u_{x}||_{L^2}^{2} + \int_0^1 \dfrac{\eta(\chi)u_{xx}^2}{v}\, dx \nonumber\\
& = \int_0^1 \left( \dfrac{\eta_x u_{x}}{v} - \dfrac{\eta u_{x}v_x}{v^2} - \dfrac{\theta_{x}}{v} + \dfrac{\theta v_x}{v} - \dfrac{1}{2}\left( \dfrac{\chi_x^2}{v^2} \right)_x \right) u_{xx}\,dx \nonumber\\
& \le \dfrac{1}{4} \int_0^1 \dfrac{\eta(\chi)u_{xx}^2}{v}\,dx + C\int_0^1 \eta^{-1}(\chi)\left( \eta_x^2 u_x^2 + \eta^2 u_x^2 v_x^2 + \theta_x^2 + \theta^2 v_x^2 + \left( \dfrac{\chi_x}{v} \right)^2 \left( \dfrac{\chi_x}{v} \right)_x^2 \right) \,dx \nonumber \\
& \le \dfrac{1}{4} \int_0^1 \dfrac{\eta(\chi)u_{xx}^2}{v}\,dx + C \int_0^1\left( \alpha^2\chi_x^2 u_x^2 + u_x^2 v_x^2 + \theta_x^2 + \theta^2 v_x^2 + \left( \dfrac{\chi_x}{v} \right)^2 \left( \dfrac{\chi_x}{v} \right)_x^2 \right) \,dx \nonumber \\
& \le \dfrac{1}{4} \int_0^1 \dfrac{\eta(\chi)u_{xx}^2}{v}\,dx + C \int_0^1\left(u_x^2 + u_x^2 v_x^2 + \theta_x^2 + \theta^2 v_x^2 + \left( \dfrac{\chi_x}{v} \right)^2 \left( \dfrac{\chi_x}{v} \right)_x^2 \right) \,dx.
\end{align}
Using $(\ref{vsxj})$, $(\ref{mu + u_x L^2})$ and $(\ref{sup v_x chi_x})$ gives
\begin{equation}\label{u_x 1}
\begin{split}
\int_0^T \int_0^1\left(u_x^2 + u_x^2 v_x^2 \right) \,dxdt & \le C + \int_0^T ||u_x||_{L^\infty}^2 ||v_x||_{L^2}^2\,dt \\
& \le C + C \int_0^T ||u_x||_{L^2} ||u_{xx}||_{L^2} \,dt \\
& \le \dfrac{1}{4} \int_0^T\int_0^1 \dfrac{\eta(\chi)u_{xx}^2}{v}\,dxdt +  C. \\
\end{split}
\end{equation}
In view of $(\ref{theta})$ and $(\ref{sup v_x chi_x})$, we derive that
\begin{align*}
\int_0^T  \int_0^1\left(\theta_x^2 + \theta^2 v_x^2 \right) \,dxdt &\le C \int_0^T  \left( ||\theta_x||_{L^2}^2 + ||\theta - \bar{\theta}+\bar{\theta}||_{L^\infty}^2 ||v_x||_{L^2}^2\right) dt \nonumber\\
& \le C \int_0^T  \left( ||\theta_x||_{L^2}^2 + ||\theta - \bar{\theta}||_{L^\infty}^2 ||v_x||_{L^2}^2 + ||v_x||_{L^2}^2 \right) dt \nonumber\\
& \le C \int_0^T ||\theta_x||_{L^2}^2 dt + C.
\end{align*}
Thanks to  $(\ref{(chi_x/v)_x L2})$ and $(\ref{sup chi_xx chi_t})$, we get
\begin{equation}\label{u_x 3}
\begin{split}
\int_0^T \int_0^1 \left( \dfrac{\chi_x}{v} \right)^2 \left( \dfrac{\chi_x}{v} \right)_x^2\,dxdt & \le \int_0^T\left|\left| \dfrac{\chi_x}{v}\right|\right|_{L^\infty}^2\left|\left| \left(\dfrac{\chi_x}{v}\right)_x\right|\right|_{L^2}^2 dt \\
& \le C \int_0^T\left|\left| \left(\dfrac{\chi_x}{v}\right)_x\right|\right|_{L^2}^2 dt \le C.\\
\end{split}
\end{equation}
Hence, integrating $(\ref{u_x 0})$ over $[0,T]$, and inserting $(\ref{u_x 1})$--$(\ref{u_x 3})$ into it, we arrive at
\begin{equation}\label{remain theta_x}
\sup \limits_{0\le t \le T}\left|\left| u_x \right|\right|_{L^2}^2 + \int_0^T \left|\left| u_{xx}  \right|\right|_{L^2}^2 dt \le C + C \int_0^T ||\theta_x||_{L^2}^2 dt.
\end{equation}

Next, we deal with the last term in $(\ref{remain theta_x})$.
Indeed, if $\beta>1$, one can choose $p=\beta-1$ in $(\ref{theta^beta theta_x})$ to obtain $||\theta_x||_{L^2(0,T;L^2)}^2 \le C$.
With which and $(\ref{remain theta_x})$, it holds that
\begin{equation}\label{u_x 11}
\sup \limits_{0\le t \le T}\left|\left| u_x \right|\right|_{L^2}^2 + \int_0^T \left(\left|\left| u_{xx}  \right|\right|_{L^2}^2 + ||\theta_x||_{L^2}^2 \right)dt \le C.
\end{equation}
In the case of $0<\beta \le 1$, multiplying $(\ref{L-NSAC1D})_5$ by $\theta$, by using $(\ref{L-NSAC1D})_1$, $(\ref{L-NSAC1D})_3$, $(\ref{manyL^1})$ and $(\ref{theta1/2_2})$, we have
\begin{align}\label{sup theta 0}
& \quad\dfrac{1}{2} \dfrac{d}{dt} ||\theta||_{L^2}^{2} + \int_0^1 \dfrac{\theta^\beta \theta_x^2}{v}\,dx \nonumber\\
& =  \int_0^1 \dfrac{\eta(\chi)u_x^2\theta}{v}\,dx -  \int_0^1 \dfrac{\theta^2 u_x}{v}\,dx +  \int_0^1 v\mu^2 \theta\,dx  \nonumber\\
& \le \int_0^1 \dfrac{\eta(\chi)u_x^2\theta}{v}\,dx -  \int_0^1 \dfrac{(\theta^2 -\bar{\theta}^2) u_x}{v}\,dx +  (1-\bar{\theta}^2)\int_0^1 \dfrac{u_x}{v}\,dx - \int_0^1 \dfrac{u_x}{v}\,dx + \int_0^1 \dfrac{\chi_t^2 \theta}{v}\,dx \nonumber\\
& \le \int_0^1 \dfrac{\eta(\chi)u_x^2\theta}{v}\,dx -  \int_0^1 \dfrac{(\theta^2 -\bar{\theta}^2) u_x}{v}\,dx + C \left| 1-\bar{\theta}^\frac{1}{2} \right| ||u_x||_{L^1} -  \left( \int_0^1 \ln v \,dx \right)_t + \int_0^1 \dfrac{\chi_t^2 \theta}{v}dx \nonumber\\
& \le \int_0^1 \dfrac{\eta(\chi)u_x^2\theta}{v}\,dx -  \int_0^1 \dfrac{(\theta^2 -\bar{\theta}^2) u_x}{v}\,dx  + \int_0^1 \dfrac{\chi_t^2 \theta}{v}\,dx - \left( \int_0^1 \ln v \,dx \right)_t + C W(t) \nonumber \\
& :=\sum\limits_{i=1}^{3}J_i  - \left( \int_0^1 \ln v \,dx \right)_t + C W(t).
\end{align}
It follows from $(\ref{v=1})$ and $(\ref{vsxj})$ that
\begin{equation}\label{sup theta 1}
\begin{split}
J_1 =  \int_0^1 \dfrac{\eta(\chi)u_x^2\theta}{v}\,dx & \le C||u_x||_{L^\infty}^2||\theta||_{L^1} \le  C||u_x||_{L^2}||u_{xx}||_{L^2}  \\
& \le \tilde{\varepsilon}_2 ||u_{xx}||_{L^2}^2 + C \tilde{\varepsilon}_2 ^{-1}||u_x||_{L^2}^2.
\end{split}
\end{equation}
For $0<\beta \le 1$, it is easy to check that $\theta^{\beta-2}+\theta^\beta \ge 1$. Thus
\begin{equation*}
\begin{split}
\left|\left| \theta^2 - \bar{\theta}^2 \right|\right|_{L^\infty} &\le C \int_0^1 |\theta \theta_x| \,dx \le C||\theta||_{L^2}||\theta_x||_{L^2} \\
& \le C ||\theta||_{L^2} \left|\left| \left( \theta^{\frac{\beta}{2}-1}+\theta^{\frac{\beta}{2}} \right) \theta_x\right|\right|_{L^2} \\
& \le C ||\theta||_{L^2} \left( W^{\frac{1}{2}}(t) +  \left|\left| \theta^{\frac{\beta}{2}} \theta_x\right|\right|_{L^2} \right).
\end{split}
\end{equation*}
Together with $(\ref{vsxj})$, $(\ref{manyL^1})$ and Cauchy-Schwartz's inequality, we get
\begin{equation}\label{sup theta 2}
\begin{split}
J_2 & = -  \int_0^1 \dfrac{(\theta^2 -\bar{\theta}^2) u_x}{v}\,dx \le C \left|\left| \theta^2 - \bar{\theta}^2 \right|\right|_{L^\infty} ||u_x||_{L^1} \\
& \le C ||\theta||_{L^2}\left( W^{\frac{1}{2}}(t) +  \left|\left| \theta^{\frac{\beta}{2}} \theta_x\right|\right|_{L^2} \right)W^{\frac{1}{2}}(t)  \\
& \le \dfrac{1}{2}\left|\left| \theta^{\frac{\beta}{2}} \theta_x\right|\right|_{L^2}^2+  C W(t)||\theta||_{L^2} + C  W(t) ||\theta||_{L^2}^2\\
& \le \dfrac{1}{2}\left|\left| \theta^{\frac{\beta}{2}} \theta_x\right|\right|_{L^2}^2+  C W(t) + C  W(t) ||\theta||_{L^2}^2.
\end{split}
\end{equation}
In view of $(\ref{v=1})$ and $(\ref{vsxj})$, we derive
\begin{equation}\label{sup theta 3}
\begin{split}
J_3 = \int_0^1 \dfrac{\chi_t^2 \theta}{v}\,dx \le C ||\chi_t||_{L^\infty}^2 \int_0^1 \theta\,dx \le C ||\chi_t||_{L^2}^2  +  C ||\chi_{xt}||_{L^2}^2.
\end{split}
\end{equation}
Hence, putting $(\ref{sup theta 1})$, $(\ref{sup theta 2})$ and $(\ref{sup theta 3})$ into $(\ref{sup theta 0})$, intergrating over $[0,T]$, and using Gronwall's inequality, we have
\begin{equation*}
\sup\limits_{0\le t \le T} ||\theta||_{L^2}^2 + \int_0^T \int_0^1 \theta^{\beta} \theta_x^2 \,dxdt \le C(\tilde{\varepsilon}_2) + \tilde{\varepsilon}_2 \int_0^T \left|\left| u_{xx}  \right|\right|_{L^2}^2 dt.
\end{equation*}
Thus
\begin{equation}\label{sup theta 4}
\begin{split}
\sup\limits_{0\le t \le T} ||\theta||_{L^2}^2 + \int_0^T \int_0^1  \theta_x^2 \,dxdt & \le \sup\limits_{0\le t  \le T} ||\theta||_{L^2}^2 + \int_0^T \int_0^1 \left( \theta^{\beta-2}+\theta^\beta \right)\theta_x^2 \,dxdt \\
& \le C(\tilde{\varepsilon}_2) + \tilde{\varepsilon}_2 \int_0^T \left|\left| u_{xx}  \right|\right|_{L^2}^2 dt.
\end{split}
\end{equation}
Plugging $(\ref{sup theta 4})$ into $(\ref{remain theta_x})$  and chossing $\tilde{\varepsilon}_2$ small enough, we conclude that $(\ref{u_x 11})$ also holds for $0<\beta \le 1$. Hence, we have shown that $(\ref{u_x 11})$ is valid for any $\beta>0$.

Finally, it follows from $(\ref{theta})$, $(\ref{sup v_x chi_x})$ and $(\ref{u_x 11})$ that
\begin{equation*}
\begin{split}
\int_0^T ||\theta v_x||_{L^2}^2\,dt & \le C \int_0^T\left( ||\theta - \bar{\theta}||_{L^\infty}^2 ||v_x||_{L^2}^2 + ||\bar{\theta}||_{L^\infty}^2||v_x||_{L^2}^2 \right)\,dt \\
& \le C \int_0^T\left( ||\theta_x||_{L^2}^2||v_x||_{L^2}^2 + ||v_x||_{L^2}^2 \right)\,dt \\
& \le C + C \int_0^T ||\theta_x||_{L^2}^2  \,dt \le C. \\
\end{split}
\end{equation*}
From which and $(\ref{vsxj})$, $(\ref{mu + u_x L^2})$, $(\ref{(chi_x/v)_x L2})$, $(\ref{sup v_x chi_x})$, $(\ref{sup chi_xx chi_t})$, $(\ref{u_x 11})$, we get
\begin{equation*}
\begin{split}
\int_0^T ||u_t||_{L^2}^2\,dt & =  \int_0^T \left|\left| \dfrac{\eta(\chi)u_{xx}}{v} + \dfrac{\eta_x u_{x}}{v} - \dfrac{\eta u_{x}v_x}{v^2} - \dfrac{\theta_{x}}{v} + \dfrac{\theta v_x}{v} - \dfrac{1}{2}\left( \dfrac{\chi_x^2}{v^2} \right)_x\right|\right|_{L^2}^2dt  \\
& \le C \int_0^T \left( ||u_{xx}||_{L^2}^2 + \alpha^2 ||\chi_x||_{L^\infty}^2||u_{x}||_{L^2}^2  + ||u_x||_{L^\infty}^2||v_{x}||_{L^2}^2 + ||\theta_{x}||_{L^2}^2 \right)\,dt \\
&\quad + C \int_0^T  \left( ||\theta v_x||_{L^2}^2 +  \left|\left| \dfrac{\chi_x}{v}\right|\right|_{L^\infty}^2\left|\left| \left(\dfrac{\chi_x}{v}\right)_x\right|\right|_{L^2}^2 \right) \,dt \\
& \le C + C \int_0^T ||\theta v_x||_{L^2}^2  \,dt \le C.
\end{split}
\end{equation*}
This, together with $(\ref{u_x 11})$, finishes the proof of  $(\ref{sup u_x 00})$.
\hfill$\Box$

\begin{Lemma}
Let the conditions of Lemma $\ref{Lmaxminv}$ be in force. Then there hold
\begin{equation}\label{max min theta}
C_1\le \theta(x,t) \le C_1^{-1},\quad \forall (x,t) \in [0,1]\times[0,T],
\end{equation}
and
\begin{equation}\label{sup theta_x L2}
\sup \limits_{0 \le t \le T} ||\theta_x(t)||_{L^2}^{2} + \int_0^T  ||(\theta_{xx},\theta_t)||_{L^2}^{2}\,dt \le C_5.
\end{equation}
\end{Lemma}

\noindent{\it\bfseries Proof.}\quad
First, we make use of $(\ref{vsxj})$ and $(\ref{sup theta 3})$ to deduce from $(\ref{sup theta 0})$ that
\begin{align}
&\quad \dfrac{1}{2} \dfrac{d}{dt} ||\theta||_{L^2}^{2} + \int_0^1 \dfrac{\theta^\beta \theta_x^2}{v}\,dx \nonumber\\
& \le \int_0^1 \dfrac{\eta(\chi)u_x^2\theta}{v}\,dx -  \int_0^1 \dfrac{(\theta^2 -\bar{\theta}^2) u_x}{v}\,dx  + \int_0^1 \dfrac{\chi_t^2 \theta}{v}\,dx - \left( \int_0^1 \ln v \,dx \right)_t + C W(t) \nonumber \\
& \le C ||u_x||_{L^\infty}^{2} +  C \left|\left| \theta^2 - \bar{\theta}^2 \right|\right|_{L^\infty} ||u_x||_{L^1} + C ||\chi_{t}||_{H^1}^2 - \left( \int_0^1 \ln v \,dx \right)_t + C W(t) \nonumber \\
& \le C ||u_x||_{H^1}^{2} +  C \left|\left| \theta^2 - \bar{\theta}^2 \right|\right|_{L^\infty} W^{\frac{1}{2}}(t) + C ||\chi_{t}||_{H^1}^2 - \left( \int_0^1 \ln v \,dx \right)_t + C W(t), \nonumber \\
& \le C ||u_x||_{H^1}^{2} +  C ||\theta||_{L^2}||\theta_x||_{L^2} W^{\frac{1}{2}}(t) + C ||\chi_{t}||_{H^1}^2 - \left( \int_0^1 \ln v \,dx \right)_t + C W(t), \nonumber \\
& \le C ||u_x||_{H^1}^{2} +  C  W(t)||\theta||_{L^2}^2 + C||\theta_x||_{L^2}^2 + C ||\chi_{t}||_{H^1}^2 - \left( \int_0^1 \ln v \,dx \right)_t + C W(t). \nonumber
\end{align}
Integrating it over $[0,T]$, using $(\ref{energy})$, $(\ref{mu + u_x L^2})$, $(\ref{sup chi_xx chi_t})$,  $(\ref{sup u_x 00})$ and Gronwall's inequality, we get
\begin{equation}\label{sup theta wxz}
\begin{split}
\sup\limits_{0\le t \le T} ||\theta||_{L^2}^2 + \int_0^T \int_0^1 \theta^{\beta} \theta_x^2 \,dxdt & \le  C.
\end{split}
\end{equation}

Next, multiplying $(\ref{L-NSAC1D})$ by $\theta^\beta \theta_t$ and integrating the result over $[0,1]$ yield
\begin{align}\label{K_0}
& \dfrac{1}{2} \dfrac{d}{dt} \int_0^1 \dfrac{(\theta^\beta \theta_x)^2}{v} dx + \int_0^1 \theta^\beta \theta_t^2\,dx \nonumber\\
& = - \dfrac{1}{2}\int_0^1 \dfrac{(\theta^\beta \theta_x)^2u_x}{v^2}dx + \int_0^1 \left(\dfrac{\eta(\chi)u_x^2-\theta u_x}{v} + v \mu^2  \right)\theta^\beta \theta_tdx  \nonumber \\
& \le \dfrac{1}{2}\int_0^1 \theta^\beta \theta_t^2\,dx + C \int_0^1 \left( \theta^{\beta+2}u_x^2 + \theta^\beta u_x^4 +\theta^\beta \mu^4 \right)dx + C ||u_x||_{L^\infty}||\theta||_{L^\infty}^{\frac{\beta}{2}}\int_0^1 \theta^{\frac{3\beta}{2}}\theta_x^2\,dx \nonumber \\
& := \dfrac{1}{2}\int_0^1 \theta^\beta \theta_t^2\,dx + K_1 + K_2.
\end{align}
The second term on the right-hand side can be estimated as follows
\begin{align}\label{K_1}
K_1 & =  C \int_0^1 \left( \theta^{\beta+2}u_x^2 + \theta^\beta u_x^4 + \theta^\beta \mu^4 \right)dx\nonumber\\
& \le C \left(||\theta||_{L^\infty}^{\beta+1}||u_x||_{L^\infty}||\theta||_{L^2}||u_x||_{L^2} +  ||\theta||_{L^\infty}^{\beta}||u_x||_{L^\infty}^2||u_x||_{L^2}^2 + ||\theta||_{L^\infty}^{\beta}||\mu||_{L^\infty}^2||\mu||_{L^2}^2 \right) \nonumber\\
& \le C\left( 1+ ||\theta||_{L^\infty}^{\beta+1} \right) ||u_{xx}||_{L^2}||u_x||_{L^2} + C \left( 1+ ||\theta||_{L^\infty}^{\beta+1} \right) ||\chi_t||_{L^\infty}^2 \\
& \le C \left( 1+ ||\theta^\beta \theta_x||_{L^2} \right) ||u_{xx}||_{L^2}||u_x||_{L^2} +  C \left( 1+ ||\theta^\beta \theta_x||_{L^2} \right) ||\chi_{t}||_{L^2}||\chi_{xt}||_{L^2} \nonumber\\
& \le C \left( ||u_x||_{L^2}^2 + ||u_{xx}||_{L^2}^2 + ||\chi_{t}||_{L^2}^2 + ||\chi_{xt}||_{L^2}^2\right) +C \left(||u_{xx}||_{L^2}^2+||\chi_{xt}||_{L^2}^2\right)||\theta^\beta \theta_x||_{L^2}^2 ,\nonumber
\end{align}
where we have used $(\ref{L-NSAC1D})_3$, $(\ref{sup chi_xx chi_t})$, $(\ref{sup u_x 00})$, $(\ref{sup theta wxz})$, Young's inequality and the following Sobolev inequalities
\begin{equation*}
||u_x||_{L^\infty} \le C ||u_{xx}||_{L^2}, \quad ||u_x||_{L^\infty}^2 \le C||u_x||_{L^2}||u_{xx}||_{L^2},
\end{equation*}
and
\begin{equation*}
||\chi_t||_{L^\infty}^2 \le C||\chi_t||_{L^2}||\chi_{xt}||_{L^2}.
\end{equation*}
Similarly, we have
\begin{equation}\label{K_2}
\begin{split}
K_2 & = C ||u_x||_{L^\infty}||\theta||_{L^\infty}^{\frac{\beta}{2}}\int_0^1 \theta^{\frac{3\beta}{2}}\theta_x^2\,dx  \\
& \le C ||u_x||_{L^\infty} \left( 1+ ||\theta||_{L^\infty}^{\beta+1} \right) ||\theta^\beta\theta_x||_{L^2}||\theta^{\frac{\beta}{2}}\theta_x||_{L^2} \\
& \le C ||u_{xx}||_{L^2}\left( 1+ ||\theta^\beta \theta_x||_{L^2}^2 \right)||\theta^{\frac{\beta}{2}}\theta_x||_{L^2} \\
& \le C \left(||u_{xx}||_{L^2}^2 + ||\theta^{\frac{\beta}{2}}\theta_x||_{L^2}^2 \right)+ C \left(||u_{xx}||_{L^2}^2 + ||\theta^{\frac{\beta}{2}}\theta_x||_{L^2}^2 \right) ||\theta^\beta\theta_x||_{L^2}^2 .
\end{split}
\end{equation}
Thus, substituting $(\ref{K_1})$ and $(\ref{K_2})$ into $(\ref{K_0})$, using $(\ref{L-NSAC1D})_3$, $(\ref{vsxj})$, $(\ref{L u_x mu L^2})$, $(\ref{sup chi_xx chi_t})$, $(\ref{sup u_x 00})$, $(\ref{sup theta wxz})$ and Gronwall's inequality, we obtain
\begin{equation}\label{sup theta_x L2}
\sup\limits_{0\le t \le T}||\theta^\beta \theta_x||_{L^2}^2 +\int_0^T \int_0^1 \theta^\beta \theta_t^2\,dxdt \le C.
\end{equation}
in view of $(\ref{theta})$, noting that
\begin{equation*}
||\theta||_{L^\infty}^{\beta+1} = ||\theta^{\beta+1}||_{L^\infty} \le ||\theta^{\beta+1}-\bar{\theta}^{\beta+1}||_{L^\infty}+ C \le C ||\theta^\beta\theta_x||_{L^2} + C \le C,
\end{equation*}
which implies
\begin{equation}\label{max theta}
\theta(x,t)\le C, \quad \forall (x,t) \in [0,1]\times[0,T].
\end{equation}

Moreover, it follows from $(\ref{sup theta wxz})$ and $(\ref{max theta})$ that
\begin{equation}\label{min theta 1}
\begin{split}
\int_0^T \int_0^1 \left( \theta^{\beta+1}-\bar{\theta}^{\beta+1} \right)^2 dxdt
&\le C \int_0^T \int_0^1 \theta^{2\beta}\theta_x^2\,dxdt \\
& \le C \sup\limits_{0\le t \le T}||\theta||_{L^\infty}^\beta \int_0^T \int_0^1 \theta^{\beta}\theta_x^2\,dxdt \le C.
\end{split}
\end{equation}
From $(\ref{energy})$, $(\ref{v=1})$, $(\ref{sup chi_xx chi_t})$ and $(\ref{sup u_x 00})$, there holds
\begin{equation*}
\begin{split}
\int_0^T \bar{\theta}_t^2 dt & = \int_0^T \left\{ \dfrac{d}{dt}\left[ 1 - \int_0^1 \left(\dfrac{u^{2}}{2}+\dfrac{(\chi^2-1)^2}{4}+\dfrac{\chi_{x}^2}{2v} \right)dx\right]\right\}^2dt \\
& = \int_0^T \left\{ \int_0^1 \left( -uu_t-(\chi^2-1)\chi\chi_t - \dfrac{\chi_{x}\chi_{xt}}{v} + \dfrac{\chi_{x}^2 v_{t}}{2v^2} \right)dx\right\}^2dt \\
& \le \int_0^T \left( ||u||_{L^2}^2||u_t||_{L^2}^2 + ||\chi||_{L^2}^2||\chi_t||_{L^2}^2 + ||\chi_x||_{L^2}^2||\chi_{xt}||_{L^2}^2 + ||\chi_x||_{L^\infty}^4||u_x||_{L^2}^2 \right)dt \\
& \le C \int_0^T \left(||u_t||_{L^2}^2 + ||\chi_t||_{L^2}^2 + ||\chi_{xt}||_{L^2}^2 + ||u_x||_{L^2}^2 \right)dt \\
& \le C.
\end{split}
\end{equation*}
Hence, we have
\begin{equation}\label{min theta 2}
\begin{split}
&\int_0^T  \left|\dfrac{d}{dt}\int_0^1 \left( \theta^{\beta+1}-\bar{\theta}^{\beta+1} \right)^2 dx\right|dt \\
\le& C \int_0^T  \int_0^1 \left( \theta^{\beta+1}-\bar{\theta}^{\beta+1} \right)^2 dxdt + C \int_0^T \left( ||\theta^\beta \theta_t||_{L^2}^{2} +  \bar{\theta}^{2\beta}\bar{\theta}_t^2 \right)dt \\
\le &C + C \int_0^T \bar{\theta}_t^2 dt \le C.
\end{split}
\end{equation}
Combining $(\ref{min theta 1})$ with $(\ref{min theta 2})$, one arrive at
\begin{equation*}
\lim \limits_{t \to +\infty} \int_0^1 \left( \theta^{\beta+1}-\bar{\theta}^{\beta+1} \right)^2 dxdt =0.
\end{equation*}
Together with $(\ref{sup theta_x L2})$, we see that, as $t\to+\infty$,
\begin{equation*}
\left|\left| \left( \theta^{\beta+1} - \bar{\theta}^{\beta+1} \right)(t) \right|\right|_{L^\infty}^2 \le C \left|\left| \left( \theta^{\beta+1} - \bar{\theta}^{\beta+1} \right)(t) \right|\right|_{L^2}^2 ||\theta^\beta \theta_x||_{L^2}^2 \to 0.
\end{equation*}
Then, by $(\ref{theta})$, we conclude that there exists a time $T_0 \gg 1$ such that
\begin{equation}\label{min theta}
\theta(x,t) \ge \dfrac{\gamma_1}{2}, \quad \forall (x,t)\in [0,1]\times [T_0,+\infty).
\end{equation}
Let $T_0$ be fixed as in $(\ref{min theta})$. Multiplying $(\ref{L-NSAC1D})_5$ by $\theta^{-p}$ with $p>2$, and integrating by parts over $[0,1]$, by $(\ref{vsxj})$, we have
\begin{equation*}
\begin{split}
& \dfrac{1}{p-1} \dfrac{d}{dt}\left|\left| \theta^{-1} \right|\right|_{L^{p-1}}^{p-1} + p \int_0^1 \dfrac{\theta^\beta \theta_x^2}{ v\theta^{p+1}} \,dx + \int_0^1 \dfrac{\eta(\chi) u_x^2}{ v\theta^{p}} \,dx + \int_0^1 \dfrac{v \mu^2}{\theta^{p}} \,dx   \\
& =  \int_0^1 \dfrac{\theta u_x}{v\theta^p} \,dx \le \dfrac{1}{2}\int_0^1 \dfrac{\eta(\chi) u_x^2}{ v\theta^{p}} \,dx + C \int_0^1 \left(\theta^{-1}\right)^{p-2}\,dx \\
& \le \dfrac{1}{2}\int_0^1 \dfrac{\eta(\chi) u_x^2}{ v\theta^{p}} \,dx + C  \left(\int_0^1 1^{p-1}\,dx\right)^{p-1} \left(\int_0^1 \left(\theta^{-1}\right)^{p-1}\,dx\right)^{\frac{p-2}{p-1}} \\
& \le \dfrac{1}{2}\int_0^1 \dfrac{\eta(\chi) u_x^2}{ v\theta^{p}} \,dx + C\left|\left| \theta^{-1} \right|\right|_{L^{p-1}}^{p-2},\\
\end{split}
\end{equation*}
which implies
\begin{equation*}
\dfrac{d}{dt}\left|\left| \theta^{-1} \right|\right|_{L^{p-1}} \le C,
\end{equation*}
where $C$ is a generic positive constant independent $p$. Hence, integrating the above inequality over $[0,t]$ and letting $p\to \infty$, we obtain
\begin{equation*}
\theta^{-1}(x,t) \le C( T_0 + 1 ) \Longleftrightarrow \theta(x,t) \ge [C( T_0 + 1 )]^{-1},\, \forall(x,t)\in [0,1]\times [0,T_0].
\end{equation*}
This, together with $(\ref{max theta})$ and $(\ref{min theta})$, proves $(\ref{max min theta})$.

Finally, using $(\ref{max min theta})$, we get from $(\ref{sup theta_x L2})$ that
\begin{equation}\label{sup theta_x L2 1}
||\theta_x||_{L^2}^2  = ||\theta^\beta \theta_x\cdot\theta^{-\beta}||_{L^2}^2 \le C_1^{-\beta}||\theta^\beta \theta_x||_{L^2}^2 \le C,
\end{equation}
and
\begin{equation}\label{sup theta_x L2 2}
\int_0^T ||\theta_{t}||_{L^2}^2 \,dt \le  C_1^{-\frac{\beta}{2}}\int_0^T||\theta^{\frac{\beta}{2}} \theta_x||_{L^2}^2  \,dt\le C.
\end{equation}
In view of $(\ref{mu + u_x L^2})$, $(\ref{sup v_x chi_x})$, $(\ref{sup chi_xx chi_t})$, $(\ref{sup u_x 00})$, $(\ref{sup theta_x L2 1})$ and $(\ref{sup theta_x L2 2})$, we deduce from $(\ref{L-NSAC1D})_5$ that
\begin{equation*}
\begin{split}
\int_0^T ||\theta_{xx}||_{L^2}^2 \,dt &\le C \int_0^T \int_0^1 \left( \theta_t^2 + u_x^2 + \theta_x^4 + \theta_x^2 v^2 + u_x^4 + \mu^4\right) \,dxdt \\
& \le C + C \int_0^T ||\theta_x||_{L^\infty}^{2} \,dt \le  C + C \int_0^T ||\theta_x||_{L^2}||\theta_{xx}||_{L^2} \,dt \\
& \le  \dfrac{1}{2}\int_0^T ||\theta_{xx}||_{L^2}^2 \,dt + C.
\end{split}
\end{equation*}
Together with $(\ref{sup theta_x L2 1})$ and $(\ref{sup theta_x L2 2})$, it leads to $(\ref{sup theta_x L2})$.
\hfill$\Box$

\begin{Lemma}
Let the conditions of Lemma $\ref{Lmaxminv}$ be in force. Then it holds that
\begin{equation}\label{sup u_t u_xx}
\sup \limits_{0 \le t \le T} ||(u_t,u_{xx})||_{L^2}^{2} + \int_0^T  ||u_{xt}||_{L^2}^{2}\,dt \le C_6.
\end{equation}
\end{Lemma}

\noindent{\it\bfseries Proof.}\quad
Differentiating $(\ref{L-NSAC1D})_2$ with respect to $t$, by $(\ref{L-NSAC1D})_1$ we find
\begin{equation*}
u_{tt} + \left( \dfrac{v\theta_t-\theta u_x}{v^2} \right)_x = \left( \left(\dfrac{\eta}{v}\right)_t u_x + \dfrac{\eta}{v}u_{xt} \right)_x  - \dfrac{1}{2}\left(\dfrac{\chi_x^2}{v^2}\right)_{xt}.
\end{equation*}
We multiply it by $u_t$ and integrate the result $[0,1]$, then
\begin{equation*}
\begin{split}
\dfrac{1}{2}& \dfrac{d}{dt}||u_t||_{L^2}^{2} + \int_0^1 \dfrac{\eta(\chi)u_{xt}^2}{v}\,dx \\
& =  - \int_0^1 \left(\dfrac{\eta}{v}\right)_t u_x u_{xt} \,dx +  \int_0^1 \dfrac{v\theta_t-\theta u_x}{v^2}u_{xt} \,dx + \int_0^1\left(\dfrac{\chi_x}{v}\right)\left(\dfrac{\chi_x}{v}\right)_{t}u_{xt}\,dx \\
& \le \dfrac{1}{2}\int_0^1 \dfrac{\eta(\chi)u_{xt}^2}{v}\,dx + C \int_0^1 \left( \chi_t^2 u_x^2 + u_x^4 + \theta_t^2 + u_x^2 + \chi_x^2\chi_{xt}^2 + \chi_x^4 u_x^2 \right)\,dx \\
& \le \dfrac{1}{2}\int_0^1 \dfrac{\eta(\chi)u_{xt}^2}{v}\,dx + C||u_x||_{L^\infty}^{2}\left( ||\chi_t||_{L^2}^{2} + ||u_x||_{L^2}^{2} + ||\chi_x||_{L^\infty}^{2}||\chi_x||_{L^2}^{2} + 1\right) \\
& \quad + C \left( ||\theta_t||_{L^2}^{2} + ||\chi_x||_{L^\infty}^{2}||\chi_{xt}||_{L^2}^{2}\right) \\
& \le \dfrac{1}{2}\int_0^1 \dfrac{\eta(\chi)u_{xt}^2}{v}\,dx + C||u_x||_{L^\infty}^{2} +  C \left( ||\theta_t||_{L^2}^{2} + ||\chi_{xt}||_{L^2}^{2}\right) \\
& \le \dfrac{1}{2}\int_0^1 \dfrac{\eta(\chi)u_{xt}^2}{v}\,dx + C\left(||u_{xx}||_{L^2}^{2} + ||\theta_t||_{L^2}^{2} + ||\chi_{xt}||_{L^2}^{2}\right),
\end{split}
\end{equation*}
where we have used $(\ref{sup chi_xx chi_t})$, $(\ref{sup u_x 00})$ and Cauchy-Schwartz's inequality. In view of $(\ref{sup chi_xx chi_t})$, $(\ref{sup u_x 00})$ and $(\ref{sup theta_x L2})$, one has
\begin{equation}\label{sup u_xx}
\sup \limits_{0 \le t \le T} ||u_{t}||_{L^2}^{2} + \int_0^T  ||u_{xt}||_{L^2}^{2}\,dt \le C.
\end{equation}
As a result, it follows from $(\ref{L-NSAC1D})_2$ that
\begin{equation*}
\begin{split}
||u_{xx}||_{L^2}^{2} & \le C \left( ||u_{t}||_{L^2}^{2} + ||\theta_{x}||_{L^2}^{2} + ||v_{x}||_{L^2}^{2} + ||u_{x}||_{L^\infty}^{2}\left(||\chi_{x}||_{L^2}^{2} + ||v_{x}||_{L^2}^{2} \right) \right) \\
& \quad + C ||\chi_{x}||_{L^\infty}^{2} \left( ||\chi_{xx}||_{L^2}^{2}+||\chi_{x}||_{L^2}^{2}||v_{x}||_{L^2}^{2} \right) \\
& \le C ||u_{x}||_{L^\infty}^{2} + C \le C ||u_x||_{L^2}||u_{xx}||_{L^2} + C \\
& \le \dfrac{1}{2}||u_{xx}||_{L^2}^{2} + C.
\end{split}
\end{equation*}
This, together with $(\ref{sup u_xx})$, finishes the proof of $(\ref{sup u_t u_xx})$.
\hfill$\Box$

\begin{Lemma}
Let the conditions of Lemma $\ref{Lmaxminv}$ be in force. Then we have
\begin{equation}\label{sup v_xx L2}
\sup \limits_{0 \le t \le T} ||v_{xx}||_{L^2}^{2} + \int_0^T  ||(v_{xx},u_{xxx})||_{L^2}^{2}\,dt \le C_7.
\end{equation}
\end{Lemma}

\noindent{\it\bfseries Proof.}\quad
Differentiating $(\ref{L-NSAC1D})_2$ with respect to $x$, we get
\begin{equation}\label{u_xt - u_xxx}
u_{xt} - \eta(\chi)\left( \dfrac{v_x}{v} \right)_{xt} = - \left( \dfrac{v\theta_x-\theta v_x}{v^2} \right)_x + \eta_x \left( \dfrac{u_x}{v} \right)_{x} + \left( \eta_x \dfrac{u_x}{v} \right)_{x} - \dfrac{1}{2}\left(\dfrac{\chi_x^2}{v^2}\right)_{xx}.
\end{equation}
Multiplying it by $\left(\dfrac{v_x}{v}\right)_x$ and integrating the result over $[0,1]$ yield
\begin{align*}
\dfrac{1}{2}& \dfrac{d}{dt}\int_0^1 \eta(\chi) \bigg| \left( \dfrac{v_x}{v} \right)_x \bigg|^2 dx + \int_0^1 \dfrac{\theta}{v}\bigg| \left( \dfrac{v_x}{v} \right)_x \bigg|^2 dx \nonumber\\
& \le C ||\chi_{t}||_{L^\infty} \int_0^1  \bigg| \left( \dfrac{v_x}{v} \right)_x \bigg|^2 dx + \dfrac{1}{2} \int_0^1 \dfrac{\theta}{v}\bigg| \left( \dfrac{v_x}{v} \right)_x \bigg|^2 dx + C ||u_{xt}||_{L^2}^2 + C \int_0^1 v_x^4 \,dx \nonumber\\
& ~~+ C \int_0^1 \left( \theta_{xx}^2 + \theta_{x}^2 v_{x}^2 + \chi_{x}^2 u_{xx}^2 + \chi_{x}^2 u_{x}^2 v_{x}^2 + \chi_{xx}^2 u_{x}^2 + \left[ \left( \dfrac{\chi_x}{v} \right)\left( \dfrac{\chi_x}{v} \right)_x \right]_x^2  \right)\,dx \nonumber \\
& \le \dfrac{1}{2} \int_0^1 \dfrac{\theta}{v}\bigg| \left( \dfrac{v_x}{v} \right)_x \bigg|^2 dx + C(||\chi_{t}||_{L^2}^2+ ||\chi_{xt}||_{L^2}^2 + ||v_{x}||_{L^2}^2 )\int_0^1 \eta(\chi) \bigg| \left( \dfrac{v_x}{v} \right)_x \bigg|^2 dx  \nonumber\\
& ~~ + C \left(||u_{xt}||_{L^2}^2 + ||v_{x}||_{L^2}^2 + ||\theta_{x}||_{H^1}^2 +||u_{xx}||_{L^2}^2 + \left|\left| \left( \dfrac{\chi_x}{v} \right)_x \right|\right|_{L^\infty}^2+ \left|\left|\left( \dfrac{\chi_x}{v} \right)_{xx} \right|\right|_{L^2}^2 \right) \nonumber\\
& \le \dfrac{1}{2} \int_0^1 \dfrac{\theta}{v}\bigg| \left( \dfrac{v_x}{v} \right)_x \bigg|^2 dx + C(||\chi_{t}||_{H^1}^2 + ||v_{x}||_{L^2}^2 )\int_0^1 \eta(\chi) \bigg| \left( \dfrac{v_x}{v} \right)_x \bigg|^2 dx \nonumber\\
&~~+ C \left( ||u_{xt}||_{L^2}^2 +||v_{x}||_{L^2}^2 + ||\theta_{x}||_{H^1}^2 +||u_{x}||_{H^1}^2 +  \left|\left| \left( \dfrac{\chi_x}{v} \right)_x \right|\right|_{L^2}^2+ \left|\left|\left( \dfrac{\chi_x}{v} \right)_{xx} \right|\right|_{L^2}^2 \right)  \nonumber\\
& \le \dfrac{1}{2} \int_0^1 \dfrac{\theta}{v}\bigg| \left( \dfrac{v_x}{v} \right)_x \bigg|^2 dx + C(||\chi_{t}||_{H^1}^2 + ||v_{x}||_{L^2}^2 )\int_0^1 \eta(\chi) \bigg| \left( \dfrac{v_x}{v} \right)_x \bigg|^2 dx  \nonumber \\
&~~+ C \left( ||u_{xt}||_{L^2}^2 + ||v_{x}||_{L^2}^2+ ||\theta_{x}||_{H^1}^2 +||u_{x}||_{H^1}^2 +  \left|\left|\chi_{x}\right|\right|_{H^1}^2+ \left|\left| \chi_t \right|\right|_{H^1}^2 \right),
\end{align*}
where we have used $(\ref{sup v_x chi_x})$, $(\ref{sup chi_xx chi_t})$ and the facts
\begin{equation*}
\begin{split}
\int_0^1 v_x^4 \,dx  & \le ||v_{x}||_{L^\infty}^2 ||v_{x}||_{L^2}^2 \le C \left|\left|\dfrac{v_{x}}{v}\right|\right|_{L^\infty}^2 ||v_{x}||_{L^2}^2  \\
& \le C \left(\left|\left|\dfrac{v_{x}}{v}\right|\right|_{L^2}^2  + \left|\left|\left(\dfrac{v_{x}}{v}\right)_x\right|\right|_{L^2}^2 \right)||v_{x}||_{L^2}^2   \\
& \le C  \left( 1 +  \int_0^1 \eta(\chi) \bigg| \left( \dfrac{v_x}{v} \right)_x \bigg|^2 dx  \right)||v_{x}||_{L^2}^2,
\end{split}
\end{equation*}
and
\begin{equation*}
\begin{split}
\left|\left|\left( \dfrac{\chi_x}{v} \right)_{xx} \right|\right|_{L^2}^2 & = \left|\left|\left( \dfrac{\chi_t}{v} +(\chi^3-\chi)\right)_{x} \right|\right|_{L^2}^2 \le C\int_0^1 \left( \chi_{xt}^2 + \chi_t^2 v_x^2 + \chi_x^2 \right)\,dx \\
& \le C \left( ||\chi_{xt}||_{L^2}^2 + ||\chi_{t}||_{L^\infty}^2 ||v_{x}||_{L^2}^2 + ||\chi_{x}||_{L^2}^2 \right) \\
& \le C \left( ||\chi_{xt}||_{L^2}^2 + ||\chi_{t}||_{L^\infty}^2  + ||\chi_{x}||_{L^2}^2 \right)\\
&\le  C \left( ||\chi_{t}||_{H^1}^2 + ||\chi_{x}||_{L^2}^2 \right).
\end{split}
\end{equation*}
Using $(\ref{mu + u_x L^2})$, $(\ref{sup v_x chi_x})$, $(\ref{sup chi_xx chi_t})$, $(\ref{sup u_x 00})$, $(\ref{sup theta_x L2})$, $(\ref{sup u_t u_xx})$ and Gronwall's inequality, we get
\begin{equation*}
\sup \limits_{0 \le t \le T} \left|\left|\left( \dfrac{v_x}{v} \right)_{x} \right|\right|_{L^2}^2 + \int_0^T  \left|\left|\left( \dfrac{v_x}{v} \right)_{x} \right|\right|_{L^2}^2 \,dt \le C.
\end{equation*}
Noting that
\begin{equation*}
\begin{split}
||v_{xx}||_{L^2}^2 & \le \left|\left|\left( \dfrac{v_x}{v} \right)_{x} v + \dfrac{v_x^2}{v} \right|\right|_{L^2}^2 \le C \left|\left|\left( \dfrac{v_x}{v} \right)_{x} \right|\right|_{L^2}^2 + C ||v_{x}||_{L^\infty}^2||v_{x}||_{L^2}^2 \\
& \le C \left|\left|\left( \dfrac{v_x}{v} \right)_{x} \right|\right|_{L^2}^2 + C ||v_{x}||_{L^\infty}^2 \le C \left|\left|\left( \dfrac{v_x}{v} \right)_{x} \right|\right|_{L^2}^2 + C ||v_{x}||_{L^2}||v_{xx}||_{L^2} \\
& \le C \left|\left|\left( \dfrac{v_x}{v} \right)_{x} \right|\right|_{L^2}^2 + \dfrac{1}{2} ||v_{xx}||_{L^2}^2 + C ||v_{x}||_{L^2}^2,\\
\end{split}
\end{equation*}
we have
\begin{equation}\label{sup v_xx}
\sup \limits_{0 \le t \le T} ||v_{xx}||_{L^2}^{2} + \int_0^T  ||v_{xx}||_{L^2}^{2}\,dt \le C.
\end{equation}
In view of  $(\ref{sup v_x chi_x})$, $(\ref{sup chi_xx chi_t})$, $(\ref{sup u_x 00})$, $(\ref{sup theta_x L2})$ $(\ref{sup u_t u_xx})$, we deduce from $(\ref{u_xt - u_xxx})$ that
\begin{equation*}
\begin{split}
\int_0^T ||u_{xxx}||_{L^2}^{2}\,dt & \le C \int_0^T \left( ||u_{xt}||_{L^2}^{2} + ||\theta_{xx}||_{L^2}^{2} + ||\theta_{x}v_x||_{L^2}^{2} + ||v_{xx}||_{L^2}^{2} + ||v_x^2||_{L^2}^{2} \right)\,dt \\
& ~~~+ C \int_0^T \left( ||\chi_x u_{xx}||_{L^2}^{2} + ||v_x u_{xx}||_{L^2}^{2}+ ||\chi_{x} u_x||_{L^2}^{2} + ||v_{x} u_x||_{L^2}^{2} \right)\,dt \\
& ~~~+ C \int_0^T \left( ||u_{xx}||_{L^2}^{2} + \left|\left|\left( \dfrac{\chi_x}{v} \right)_{x}^2 \right|\right|_{L^2}^2  + \left|\left| \dfrac{\chi_x}{v}\right|\right|_{L^\infty}^2\left|\left|\left( \dfrac{\chi_x}{v} \right)_{xx} \right|\right|_{L^2}^2   \right)\,dt\le C.
\end{split}
\end{equation*}
It combined with $(\ref{sup v_xx})$ lead to $(\ref{sup v_xx L2})$.
\hfill$\Box$

\begin{Lemma}\label{sup chi_xxx}
Let the conditions of Lemma $\ref{Lmaxminv}$ be in force. Then we have
\begin{equation}\label{sup chi_xt chi_xxx}
\sup \limits_{0 \le t \le T} \left|\left|\left(\chi_{xxx},\chi_{xt},\left( \dfrac{\chi_x}{v} \right)_t\, \right)\right|\right|_{L^2}^{2} + \int_0^T \left|\left|\left(\chi_{xxx},\chi_{xxt},\left( \dfrac{\chi_t}{v} \right)_t,\left( \dfrac{\chi_x}{v} \right)_{xt}\, \right)\right|\right|_{L^2}^{2} dt \le C_8.
\end{equation}
\end{Lemma}

\noindent{\it\bfseries Proof.}\quad
First, differentiating $(\ref{chi_t/v})$ with respect to $t$, one has
\begin{equation}\label{(chi_t/v)_t}
\left(\dfrac{\chi_t}{v}\right)_t =  \left( \dfrac{\chi_x}{v} \right)_{tx} - (3\chi^2-1)\chi_t.
\end{equation}
Multiplying it by $\left(\frac{\chi_t}{v}\right)_t$ and integrating the result over $[0,1]$ yield
\begin{equation}\label{(chi_t/v)_t 1}
\begin{split}
\left|\left| \left(\dfrac{\chi_t}{v}\right)_t \right|\right|_{L^2}^2  & = - \int_0^1 \left( \dfrac{\chi_x}{v} \right)_{t} \left(\dfrac{\chi_t}{v}\right)_{xt}\,dx  - \int_0^1 (3\chi^2-1)\chi_t \left(\dfrac{\chi_t}{v}\right)_t \,dx \\
& = -  \int_0^1 \left( \dfrac{\chi_x}{v} \right)_{t} \left(\dfrac{\chi_x}{v}\right)_{tt}\,dx
-  \int_0^1 \left( \dfrac{\chi_x}{v} \right)_{t} \left(\dfrac{\chi_{x}v_t}{v^2} - \dfrac{\chi_{t}v_x}{v^2} \right)_t \,dx\\
& \quad - \int_0^1 (3\chi^2-1)\chi_t \left(\dfrac{\chi_t}{v}\right)_t \,dx,
\end{split}
\end{equation}
where we have used the facts
\begin{equation*}
\begin{split}
\left(\dfrac{\chi_t}{v}\right)_x = \left(\dfrac{\chi_x}{v}\right)_t + \dfrac{\chi_{x}v_t}{v^2} - \dfrac{\chi_{t}v_x}{v^2} \quad {\rm and}\quad
\left(\dfrac{\chi_x}{v}\right)_t = \left(\dfrac{\chi_t}{v}\right)_x + \dfrac{\chi_{t}v_x}{v^2} - \dfrac{\chi_{x}v_t}{v^2}.
\end{split}
\end{equation*}
Thus, it follows from $(\ref{(chi_t/v)_t 1})$ that
\begin{align}\label{chi_xt chi_tt 0}
\dfrac{1}{2}&\dfrac{d}{dt}\left|\left| \left(\dfrac{\chi_x}{v}\right)_t \right|\right|_{L^2}^2 +
\left|\left| \left(\dfrac{\chi_t}{v}\right)_t \right|\right|_{L^2}^2  \nonumber\\
& = -  \int_0^1 \left( \dfrac{\chi_x}{v} \right)_{t} \left(\dfrac{\chi_{x}v_t}{v^2} - \dfrac{\chi_{t}v_x}{v^2} \right)_t \,dx - \int_0^1 (3\chi^2-1)\chi_t \left(\dfrac{\chi_t}{v}\right)_t \,dx \nonumber\\
&  = -  \int_0^1 \left( \dfrac{\chi_x}{v} \right)_{t} \left( \dfrac{\chi_{xt}v_t + \chi_xv_{tt}}{v^2}-\dfrac{2\chi_{x}v_t v_x}{v^3} - \left(\dfrac{\chi_{t}}{v}\right)_t \dfrac{v_{x}}{v}-\left(\dfrac{\chi_{t}}{v}\right )\dfrac{v_{xt}}{v} + \left(\dfrac{\chi_{t}}{v}\right) \dfrac{v_{x}v_{t}}{v^2}\right)\,dx \nonumber \\
& \quad - \int_0^1 (3\chi^2-1)\chi_t \left(\dfrac{\chi_t}{v}\right)_t \,dx \nonumber\\
& \le \dfrac{1}{4}\left|\left| \left(\dfrac{\chi_t}{v}\right)_t \right|\right|_{L^2}^2 + C\left|\left| v_x\left(\dfrac{\chi_x}{v}\right)_t \right|\right|_{L^2}^2  + C \left( ||u_x||_{L^\infty}^{2} + ||\chi_x||_{L^\infty}^{2} + ||\chi_t||_{L^\infty}^{2} \right) \left|\left| \left(\dfrac{\chi_x}{v}\right)_t \right|\right|_{L^2}^2 \nonumber \\[0.3em]
& \quad + C \left( ||\chi_t||_{L^2}^{2} + ||\chi_{xt}||_{L^2}^{2} + ||u_{xt}||_{L^2}^{2} + ||u_x v_x||_{L^2}^{2} + ||u_{xx}||_{L^2}^{2} \right)  \nonumber \\[0.3em]
& \le \dfrac{1}{4}\left|\left| \left(\dfrac{\chi_t}{v}\right)_t \right|\right|_{L^2}^2  + C \left|\left| \left(\dfrac{\chi_x}{v}\right)_t \right|\right|_{L^\infty}^2 + C \left( ||u_x||_{H^1}^{2} + ||\chi_{xx}||_{L^2}^{2} + ||\chi_t||_{H^1}^{2} \right) \left|\left| \left(\dfrac{\chi_x}{v}\right)_t \right|\right|_{L^2}^2 \nonumber \\[0.3em]
& \quad + C \left( ||\chi_t||_{L^2}^{2} + ||\chi_{xt}||_{L^2}^{2} + ||u_{xt}||_{L^2}^{2} + ||u_{x}||_{H^1}^{2} \right),
\end{align}
where we have used $(\ref{sup v_x chi_x})$ and Cauchy-Schwartz's inequality. The second term on the right-hand side can be estimated as follows
\begin{equation*}
\begin{split}
\left|\left| \left(\dfrac{\chi_x}{v}\right)_t \right|\right|_{L^\infty}^2  & \le C \left|\left| \left(\dfrac{\chi_x}{v}\right)_t \right|\right|_{L^2} \left|\left| \left(\dfrac{\chi_x}{v}\right)_{xt} \right|\right|_{L^2} \\
& = C \left|\left| \left(\dfrac{\chi_x}{v}\right)_t \right|\right|_{L^2} \left|\left| \left(\dfrac{\chi_t}{v}+ (\chi^3-\chi)\right)_{t} \right|\right|_{L^2} \\
& \le \dfrac{1}{4}\left|\left| \left(\dfrac{\chi_t}{v}\right)_t \right|\right|_{L^2}^2 + C \left|\left| \left(\dfrac{\chi_x}{v}\right)_t \right|\right|_{L^2}^2 + C ||\chi_t||_{L^2}^{2}.
\end{split}
\end{equation*}
Putting it into $(\ref{chi_xt chi_tt 0})$, using $(\ref{mu + u_x L^2})$, $(\ref{sup v_x chi_x})$, $(\ref{sup chi_xx chi_t})$, $(\ref{sup u_x 00})$, $(\ref{sup u_t u_xx})$ and Gronwall's inequality yield
\begin{equation}\label{sup chi_xt chi_xxx 1}
\sup \limits_{0 \le t \le T} \left|\left| \left( \dfrac{\chi_x}{v} \right)_t\right|\right|_{L^2}^{2} + \int_0^T \left|\left|\left(\dfrac{\chi_t}{v} \right)_t \right|\right|_{L^2}^{2} dt \le C,
\end{equation}
which implies
\begin{equation*}
||\chi_{xt}||_{L^2}^{2} \le \left|\left| \left( \dfrac{\chi_x}{v} \right)_t v + \dfrac{\chi_x u_x}{v}\right|\right|_{L^2}^2 \le C \left|\left| \left( \dfrac{\chi_x}{v} \right)_t\right|\right|_{L^2}^{2} + C ||\chi_{x}||_{L^\infty}^{2}||u_{x}||_{L^2}^{2} \le C.
\end{equation*}
Moreover, in view of $(\ref{sup v_xx L2})$, we obtain
\begin{align*}
||\chi_{xxx}||_{L^2}^{2} & = \left|\left|\chi_{xt} + \dfrac{\chi_{xx}v_x}{v}+ \dfrac{\chi_{x}v_{xx}}{v} - \dfrac{\chi_{x}v_x^2}{v^2}+ v_x(\chi^3-\chi)+ v(3\chi^2-1)\chi_x \right|\right|_{L^2}^{2} \nonumber \\
& \le C \left( ||\chi_{xt}||_{L^2}^{2} + ||v_{x}||_{L^\infty}^{2} + ||v_{xx}||_{L^2}^{2} + ||v_{x}||_{L^2}^{2}+||\chi_{x}||_{L^2}^{2} \right) \nonumber \\
& \le C \left( ||\chi_{xt}||_{L^2}^{2} + ||v_{x}||_{H^1}^{2} + ||\chi_{x}||_{L^2}^{2} \right) \le C.
\end{align*}
From which we have
\begin{align*}
\int_0^T ||\chi_{xxx}||_{L^2}^{2}\,dt \le C \int_0^T \left( ||\chi_{xt}||_{L^2}^{2} + ||v_{x}||_{H^1}^{2} + ||\chi_{x}||_{L^2}^{2} \right) \,dt \le C.
\end{align*}
Finally, it follows from $(\ref{(chi_t/v)_t})$ that
\begin{align}\label{sup chi_xt chi_xxx 5}
\int_0^T \left|\left| \left( \dfrac{\chi_x}{v} \right)_{xt}\right|\right|_{L^2}^2 dt \le C \int_0^T \left( \left|\left|\left( \dfrac{\chi_t}{v} \right)_{t}\right|\right|_{L^2}^2 + ||\chi_{t}||_{L^2}^{2} \right)dt \le C,
\end{align}
which gives
\begin{equation*}
\begin{split}
\int_0^T & ||\chi_{xxt}||_{L^2}^{2}\,dt  = \int_0^T \left|\left| \left( \dfrac{\chi_x}{v} \right)_{xt}v + \dfrac{\chi_{xx}v_t}{v} + \dfrac{\chi_{xt}v_x}{v} + \dfrac{\chi_{x}v_{xt}}{v} - \dfrac{2\chi_{x}v_x v_t}{v^2}\right|\right|_{L^2}^{2} \,dt \\
& \le C \int_0^T \left( \left|\left| \left( \dfrac{\chi_x}{v} \right)_{xt}\right|\right|_{L^2}^2 + ||u_{x}||_{L^\infty}^{2} + ||\chi_{xt}||_{L^\infty}^{2} + ||\chi_{x}||_{L^\infty}^{2} + ||\chi_{x}u_x||_{L^\infty}^{2}\right)dt \\
& \le C + C \int_0^T ||\chi_{xt}||_{L^2}||\chi_{xxt}||_{L^2} \,dt \le C + \dfrac{1}{2} \int_0^T  ||\chi_{xxt}||_{L^2}^{2}\,dt.
\end{split}
\end{equation*}
This, together with $(\ref{sup chi_xt chi_xxx 1})$--$(\ref{sup chi_xt chi_xxx 5})$, leads to $(\ref{sup chi_xt chi_xxx})$.
\hfill$\Box$

\begin{Lemma}
Let the conditions of Lemma $\ref{Lmaxminv}$ be in force. Then it holds that
\begin{equation}\label{sup theta_xx theta_t}
\sup \limits_{0 \le t \le T} \left|\left|\left(\theta_{xx},\theta_{t}\right)\right|\right|_{L^2}^{2} + \int_0^T \left|\left|\left(\theta_{xxx},\theta_{xt}\right)\right|\right|_{L^2}^{2} dt \le C_9.
\end{equation}
\end{Lemma}

\noindent{\it\bfseries Proof.}\quad
First, differentiating $(\ref{L-NSAC1D})_5$ with respect to $t$, by $(\ref{L-NSAC1D})_1$ we get
\begin{equation*}
\theta_{tt} - \left( \dfrac{\theta^\beta \theta_{xt}}{v} + \dfrac{\beta\theta^{\beta-1}\theta_t \theta_x}{v} - \dfrac{\theta^{\beta}\theta_x u_x}{v^2} \right)_{x} = - \left( \dfrac{\theta u_x}{v}\right)_t  + \left( \dfrac{\eta(\chi)u_x^2}{v}\right)_t + \left( v\mu^2 \right)_t.
\end{equation*}
Multiplying it by $\theta_t$ and integrating the result over $[0,1]$ yield
\begin{align*}
\dfrac{1}{2}& \dfrac{d}{dt} ||\theta_{t}||_{L^2}^{2} + \int_0^1 \dfrac{\theta^\beta \theta_{xt}^2}{v} \,dx \nonumber\\
& \le \dfrac{1}{2}\int_0^1 \dfrac{\theta^\beta \theta_{xt}^2}{v} \,dx + C \int_0^1 \left( \theta_t^2\theta_x^2 + \theta_x^2 u_x^2 \right)dx + C  \int_0^1 \left( \theta_t^2|u_x| + |u_{xt}\theta_t| + u_x^2|\theta_t| \right)dx \nonumber\\
& \quad + C \int_0^1 \left( |\chi_t|u_x^2 + |u_x u_{xt}| + |u_x|^3 + |u_x|\mu^2 + |\mu\mu_t| \right)|\theta_t|\,dx \nonumber \\
& \le \dfrac{1}{2}\int_0^1 \dfrac{\theta^\beta \theta_{xt}^2}{v} \,dx + C \left( ||\theta_{x}||_{L^\infty}^{2}||\theta_{t}||_{L^2}^{2} + ||u_{x}||_{L^\infty}^2 ||\theta_{x}||_{L^2}^{2} + ||u_{x}||_{L^\infty}||\theta_{t}||_{L^2}^{2} \right) \nonumber \\
& \quad + C \left( ||u_{xt}||_{L^2}^{2}+ ||u_{x}||_{L^2}^{2} + ||\chi_{t}u_x||_{L^2}^{2} + ||u_{x}^2||_{L^2}^{2} + ||\mu^2||_{L^2}^{2} + ||\mu_t||_{L^2}^{2} \right)  \nonumber \\
& \quad + C\left( 1 + ||u_{x}||_{L^\infty}^2 + ||\mu||_{L^\infty}^2 \right) ||\theta_{t}||_{L^2}^2 \nonumber \\
& \le \dfrac{1}{2}\int_0^1 \dfrac{\theta^\beta \theta_{xt}^2}{v} \,dx  +  C \left( 1 + ||\theta_{x}||_{L^\infty}^2 +||u_{x}||_{L^\infty}^2 + ||\mu||_{L^\infty}^2 \right) ||\theta_{t}||_{L^2}^2  \nonumber \\
& \quad + C \left( ||u_{x}||_{H^1}^2 + ||u_{xt}||_{L^2}^{2} + ||\chi_{t}||_{L^\infty}^{2} + ||\mu||_{L^\infty}^{2} + ||\mu_t||_{L^2}^{2} \right) \nonumber \\
& \le \dfrac{1}{2}\int_0^1 \dfrac{\theta^\beta \theta_{xt}^2}{v} \,dx +  C \left( ||\theta_{x}||_{H^1}^2 +||u_{x}||_{H^1}^2 + ||\chi_t||_{H^1}^2 \right) ||\theta_{t}||_{L^2}^2 \nonumber\\
& \quad +  C \left( ||\theta_{t}||_{L^2}^2 + ||u_{x}||_{H^1}^2 + ||u_{xt}||_{L^2}^{2} + ||\chi_{t}||_{H^1}^{2} + \left|\left|\left(\dfrac{\chi_t}{v}\right)_t\right|\right|_{L^2}^{2} \right),
\end{align*}
where we have used $(\ref{sup chi_xx chi_t})$, $(\ref{sup u_x 00})$ and $(\ref{sup theta_x L2})$.
Combining the above inequality with $(\ref{mu + u_x L^2})$, $(\ref{sup chi_xx chi_t})$, $(\ref{sup u_x 00})$, $(\ref{sup theta_x L2})$, $(\ref{sup u_t u_xx})$, $(\ref{sup chi_xt chi_xxx})$ and Gronwall's inequality, we have
\begin{equation}\label{sup theta_xx theta_t 1}
\sup \limits_{0 \le t \le T} ||\theta_{t}||_{L^2}^{2} + \int_0^T ||\theta_{xt}||_{L^2}^{2}\,dt \le C.
\end{equation}
Next, it follows from $(\ref{L-NSAC1D})_5$ that
\begin{equation}\label{sup theta_xx theta_t 2}
\begin{split}
||\theta_{xx}||_{L^2}^{2} & \le C  \int_0^1 \left( \theta_x^4 + \theta_x^2 v_x^2 + \theta_t^2 + u_x^2 + u_x^4 + \mu^4 \right)dx \\
& \le C  \left( ||\theta_{x}||_{L^\infty}^{2} + ||u_{x}||_{L^\infty}^{2} + ||\chi_t||_{L^\infty}^{2} \right) + C \\
& \le C ||\theta_{x}||_{L^2}||\theta_{xx}||_{L^2} + C \le \dfrac{1}{2}||\theta_{xx}||_{L^2}^{2} + C.
\end{split}
\end{equation}
Finally, differentiating $(\ref{L-NSAC1D})_5$ with respect to $x$, by $(\ref{L-NSAC1D})_1$ we obtain
\begin{equation*}
\dfrac{\theta^\beta \theta_{xxx}}{v} =  - \theta_{xx} \left( \dfrac{\theta^\beta }{v} \right)_x + \left( -\dfrac{\beta \theta^{\beta-1}\theta_x^2}{v} + \dfrac{\theta^{\beta}\theta_x v_x}{v^2} + \theta_t + \dfrac{\theta}{v}u_x - \dfrac{\eta(\chi)u_x^2}{v} - v\mu^2 \right)_x.
\end{equation*}
Thus, there holds
\begin{equation*}
\begin{split}
\int_0^T ||\theta_{xxx}||_{L^2}^2dt & \le C \int_0^T \int_0^1 \left( \theta_x^2\theta_{xx}^2 + v_x^2\theta_{xx}^2 + \theta_x^6 + \theta_x^4 v_{x}^2 + \theta_x^2 v_{xx}^2 + \theta_x^2 v_{x}^4\right) dxdt \\
& \quad + C \int_0^T \int_0^1 \left( \theta_{xt}^2 +  \theta_x^2 u_{x}^2 +  u_{xx}^2 + u_x^2 v_{x}^2 + \chi_x^2 u_x^4 + u_x^2 u_{xx}^2 \right) dxdt \\
& \quad + C \int_0^T \int_0^1 \left( u_x^4 v_{x}^2 + v_x^2\mu^4 + \mu^2 \mu_x^2 \right) dxdt \\
& \le C.
\end{split}
\end{equation*}
From which and $(\ref{sup theta_xx theta_t 1})$, $(\ref{sup theta_xx theta_t 2})$, we obtain $(\ref{sup theta_xx theta_t})$.
\hfill$\Box$

\section{Proof of Theorem $\ref{th:1.1}$}\label{section proof}
\setcounter{equation}{0}

First, the local existence of unique strong solutions can be proved by Banach Fixed Point theorem on (see \cite{Tani}).
\begin{Lemma}\label{lemma 3.1}
Suppose that {\rm(}$\ref{local C}${\rm)} holds. Then there exists $T_0=T_0(V_0,V_0,V_0,M_0)$, depending only on $\beta,V_0$ and $M_0$, such that the initial boundary value problem $(\ref{L-NSAC1D})$--$(\ref{BC})$ has a unique solution
\begin{equation*}
(v,u,\chi,\theta)\in X(0,T_0;\frac{1}{2}V_0,V_0,\frac{1}{2}V_0,KM_0),
\end{equation*}
where $K\ge2$ is a positive constant.
\end{Lemma}

In what follows, by the a priori estimates established in Section $\ref{section priori}$ and local well-posedness result,
we derive the existence and uniqueness of global strong solutions.

\noindent{\it\bfseries Proof of Theorem $\ref{th:1.1}$.}\quad
It follows from Lemma $\ref{lemma 3.1}$ that the initial boundary value problem $(\ref{L-NSAC1D})$--$(\ref{BC})$ has a unique solution
$$
(v,u,\chi,\theta) \in X(0,t_1;\frac{1}{2}V_0,V_0,\frac{1}{2}V_0,KM_0),
$$
where $t_1=T_0(V_0,V_0,V_0, M_0)>0$.

Moreover, if we take $\alpha\le\alpha_1$ with $\alpha_1$ being small enough such that
\begin{align}
\label{local 0}
{\Big(V_0\Big)^{-\alpha_1}\le2}, \quad
{(2KM_0)^{\alpha_1}\le 1}, \quad
{\alpha_1 H\Big(\frac12 V_0, V_0, \frac12 V_0, KM_0\Big)\le\varepsilon_1},
\end{align}
where $\varepsilon_1>0$ is chosen in Lemma $\ref{Lmaxminv}$.
Then we deduce from Lemma $\ref{L-energy}$--$\ref{sup chi_xxx}$ that the solution $(v,u,\chi,\theta)$ satisfies
\begin{equation}\label{local 1}
C_0\le v\le C_0^{-1}, \quad
V_0\le \chi \le 1, \quad
C_1\le \theta \le C_1^{-1}, \quad
(x,t)\in [0,1]\times[0,t_1],
\end{equation}
and
\begin{equation}\label{local 2}
\sup\limits_{0\le t \le t_1} ||(v,u,\chi, \theta)(t)||_{H^2}^2
+ \int_0^{t_1}||\chi_t||_{L^2}^2 dt \le C_{10}^2 := \sum_{i=2}^{9} C_i,
\end{equation}
where $C_2,\cdots,C_9$ are the same ones as in Section $\ref{section priori}$.

Next, taking $(v,u,\chi,\theta)(\cdot,t_1)$ as the initial data and applying Lemma $\ref{lemma 3.1}$ again, we can extend the local solution $(v,u,\chi,\theta)$ to the time interval $[t_1,t_1+t_2]$ with $t_2=T_0(C_0,V_0,C_1$, $C_{10})$. Moreover, we have
\begin{equation}\label{local 3}
v\ge \dfrac{1}{2}C_0,\quad  \chi \ge V_0,\quad \theta \ge \dfrac{1}{2}C_0,
\qquad (x,t)\in [0,1]\times[t_1,t_1+t_2],
\end{equation}
and
\begin{equation}\label{local 4}
\sup\limits_{t_1\le t \le t_1+t_2} ||(v,u,\chi, \theta)(t)||_{H^2}^2
+ \int_{t_1}^{t_1+t_2}||\chi_t||_{L^2}^2 dt \le (KC_{10})^2,
\end{equation}
Hence, collecting ($\ref{local 1}$)--($\ref{local 4}$), we get
\begin{equation*}
v\ge \dfrac{1}{2}C_0,\quad  \chi \ge V_0,\quad \theta \ge \dfrac{1}{2}C_0,
\qquad (x,t)\in [0,1]\times[0,t_1+t_2],
\end{equation*}
and
\begin{equation*}
\sup\limits_{0 \le t \le t_1+t_2}||(v,u,\chi, \theta)(t)||_{H^2}^2
+ \int_{0}^{t_1+t_2}||\chi_t||_{L^2}^2 dt \le (K^2+1)C_{10}^2.
\end{equation*}

Take $\alpha \le \min\{\alpha_1,\alpha_2\}$, where $\alpha_1$ is the same as one as in ($\ref{local 0}$) and $\alpha_2$ is chosen to be such that
\begin{equation*}
\Big(V_0\Big)^{-\alpha_2}\le2, \quad
\Big(2\sqrt{1+K^2}C_{10}\Big)^{\alpha_2}\le 1, \quad
\alpha_2 H\Big(\frac12C_0,\, V_0,\, \frac12 C_1, \,\sqrt{1+K^2}C_{10}\Big)\le\varepsilon_1,
\end{equation*}
where $\varepsilon_1>0$ is chosen in Lemma $\ref{Lmaxminv}$. Then we infer from Lemma $\ref{L-energy}$--$\ref{sup chi_xxx}$ again that the solution $(v,u,\chi,\theta)$ satisfies ($\ref{local 1}$) and ($\ref{local 2}$) on $[0,t_1+t_2]$.

Hence, chossing $\epsilon_0 \triangleq \min\{\alpha_1,\alpha_2\}$ and repeating the above procedure, we see that the initial boundary value problem $(\ref{L-NSAC1D})$--$(\ref{BC})$ has a unique solution $(v,u,\chi,\theta)\in X(0,+\infty;C_0,V_0,C_1,C_{10})$.
This completes the proof of the global existence of strong solutions.
The uniqueness of the solutions can be easily obtained by the standard energy method.
\hfill$\Box$

\section*{Acknowledgments}

Ding's work is supported by the Key Project of National Natural Science Foundation of China (No. 12131010) the National Natural Science Foundation of China (No. 11571117, 11871005, 11771155), the Guangdong Provincial Natural Science Foundation (No. 2017A030313003, 2021A1515010249,  2021A1515010303) and the Science and Technology Program of Guang-\\zhou (No. 2019050001).
Li's work is supported by the National Natural Science Foundation of China (No.11471127, 11671155, 11771156), the Guangdong Provincial Natural Science Foundation (No. 2016A030313418) and the Science and Technology Program of Guangzhou (No. 201607010207, 201707010136).

\appendix
\section{Appendix}
\setcounter{equation}{0}
In the proof of local existence, it is easy to show that $v \in L^{\infty}(Q_T)$, $A^{-1}\le v \le A$
(for some constant $A>1$ and $Q_T=(0,t)\times(0,T)$), which leads to the following lemma.

\begin{Lemma}\label{A-chi<1}
Assume $v \in L^{\infty}(Q_T)$, $A^{-1}\le v \le A$, $\chi$ is a smooth solution to $(\ref{L-NSAC1D})_3$--$(\ref{L-NSAC1D})_4$ with $V_0 \le \chi_0 \le 1$ and $\chi_{x}\big|_{x=0,1}=0$. Then
\begin{equation*}
V_0 \le \chi(x,t) \le 1,\quad (x,t)\in [0,1]\times[0,T_0].
\end{equation*}
\end{Lemma}

\noindent{\it\bfseries Proof.}\quad
We rewrite the equations $(\ref{L-NSAC1D})_3$ and $(\ref{L-NSAC1D})_4$ as
\begin{equation*}
\chi_t = \chi_{xx} - \dfrac{v_x}{v}\chi_x - v \chi( \chi^2 - 1 ).
\end{equation*}
Multiplying the above equation by $2\chi$ gives
\begin{equation*}
(\chi^2 - 1)_{t} - (\chi^2 - 1)_{xx} + \dfrac{v_x}{v}(\chi^2 - 1)_{x} + 2v(\chi^2 - 1)= - 2v(\chi^2 - 1)^2 - 2\chi_x^2 \le 0.
\end{equation*}
Then the maximum principle implies $\chi^2 - 1 \le 0$. Thus, we have
\begin{equation*}
\chi(x,t)\le 1.
\end{equation*}

Define $Y(x,t)= e^{At}\chi(x,t)$. Then $Y$ satisfies
\begin{equation*}
Y_t - Y_{xx} +  \dfrac{v_x}{v}Y_{x} + [A + v(\chi^2-1)]Y = 0.
\end{equation*}
By the maximum principle and the definition of $Y$, we arrive at
\begin{equation*}
\chi(x,t)\ge 0.
\end{equation*}

Then consider the equation
\begin{equation*}
(\chi - V_0)_t - (\chi - V_0)_{xx} + \dfrac{v_x}{v}(\chi - V_0)_x  = - v \chi(\chi^2 - 1) \ge 0.
\end{equation*}
The maximum principle implies that
\begin{equation*}
\chi(x,t)\ge V_0.
\end{equation*}
This completes the proof of Lemma $\rm{\ref{A-chi<1}}$.
\hfill$\Box$


\end{document}